\documentclass[12pt]{article}
\usepackage{amsfonts}
\usepackage{mathrsfs}
\usepackage{amsthm}
\usepackage{amsmath}
\usepackage{amssymb}
\usepackage{latexsym}
\usepackage{graphicx}
\usepackage{subfigure}

\usepackage{color,varioref}
\usepackage{longtable}
\usepackage{rotating,multirow,array}
\usepackage{enumerate}
\usepackage{caption}
\usepackage{float}
\definecolor{red}{rgb}{0.9,0,0}

\definecolor{green}{rgb}{0,0.9,0}

\definecolor{blue}{rgb}{0,0,0.9}

\def\mc{\multicolumn}

 \def\R{{\cal R}}

\def\norm#1{\|#1\|}

\def\ni{\noindent}

\newtheorem{theorem}{Theorem}[section]
\newtheorem{proposition}{Proposition}[section]
\newtheorem{lemma}{Lemma}[section]

\newtheorem{remark}{Remark}[section]
\newtheorem{definition}{Definition}[section]
\newtheorem{assumption}{Assumption}[section]

\begin{document}

\title{\bf A proximal-proximal majorization-minimization algorithm for nonconvex tuning-free robust regression problems}

\author{
Peipei Tang\thanks{School of Computer and
Computing Science, Zhejiang University City
College, Hangzhou 310015, China (Email: tangpp@zucc.edu.cn). This author's research is supported by the Natural Science Foundation of Zhejiang Province of China under Grant No. LY19A010028, the Zhejiang Science and Technology Plan Project of China (No. 2020C03091, No. 2021C01164) and the Scientific Research Foundation of Zhejiang University City College (No. X-202112).},
Chengjing Wang
\thanks{{\bf Corresponding author}, School of Mathematics, Southwest Jiaotong University, No.999, Xian Road, West Park, High-tech Zone, Chengdu 611756, China ({Email: renascencewang@hotmail.com}). } \ and Bo Jiang\thanks{School of Computer Science, Zhejiang University, No. 38, Zheda Road, Hangzhou 310027, China ({Email: 22021105@zju.edu.cn}). }
}

\maketitle


\begin{abstract}
In this paper, we introduce a proximal-proximal majorization-minimization (PPMM) algorithm for nonconvex tuning-free robust regression problems. The basic idea is to apply the proximal majorization-minimization algorithm to solve the nonconvex problem with the inner subproblems solved by a sparse semismooth Newton (SSN) method based proximal point algorithm (PPA). We must emphasize that the main difficulty in the design of the algorithm lies in how to overcome the singular difficulty of the inner subproblem. Furthermore, we also prove that the PPMM algorithm converges to a d-stationary point. Due to the Kurdyka-{\L}ojasiewicz (KL) property of the problem, we present the convergence rate of the PPMM algorithm. Numerical experiments demonstrate that our proposed algorithm outperforms the existing state-of-the-art algorithms.
\end{abstract}


\section{Introduction}
\label{sec:Introduction}
During the last decade, the estimation of high dimensional sparse statistical models has attracted tremendous interests, in which the number of variables is probably larger than the sample size. Consider a high-dimensional linear regression model
\begin{eqnarray*}
b=X\dot{\beta}+\varepsilon,
\end{eqnarray*}
where $b\in\mathcal{R}^{n}$ is an $n$-dimensional response vector, $X\in\mathcal{R}^{n\times p}$ is a design matrix, $\varepsilon\in\mathcal{R}^{n}$ is a random error vector. Sparse learning based methods, which aim to minimize the fitting errors along with some sparse regularization terms, have received considerable attention in recent years due to their good performance and interpretability. The $\ell_{1}$ regularized least square regression that is called LASSO \cite{Tibshirani1996} requires solving the following optimization problem
\begin{eqnarray*}
\min_{\beta\in\mathcal{R}^{p}}\Big\{\frac{1}{2n}\|X\beta-b\|^{2}+\lambda\|\beta\|_{1}\Big\},
\end{eqnarray*}
where $\|\cdot\|$ is the Euclidean norm in $\R^n$ and $\lambda$ denotes the tuning parameter which controls the complexity of the model. The Lasso estimator can achieve a near-oracle performance if the error vector $\varepsilon$ is normally distributed and some suitable design conditions hold. The theory of the Lasso reveals that the tuning parameter $\lambda$ relies on the deviation of the noise in order to recover $\dot{\beta}$. However, it is of great challenge to estimate the deviation for large-scale problems, especially when the feature dimension $p$ is much larger than the sample size $n$. The square-root Lasso estimator
\begin{eqnarray*}
\min_{\beta\in\R^{p}}\Big\{\norm{X
\beta-b}+\lambda\norm{\beta}_{1}\Big\}
\end{eqnarray*}
was proposed in \cite{BelloniCW2011} to avoid the aforementioned challenge. It is known (see e.g., \cite{Bellec,Derumigny}) that the square-root Lasso estimator can achieve the minimax optimal rate of convergence under some suitable conditions which are independent of the unknown noise level. Another obstacle is that the Gaussian or sub-Gaussian error assumption can hardly be satisfied for high-dimensional microarray data, climate data, insurance claim data, e-commerce data and many other applications due to the heavy-tailed errors, which can affect the choice of $\lambda$ and result in misleading results if we apply the standard procedures directly. Therefore, Wang et al. \cite{WangPeng} studied the following $\ell_{1}$ regularized tuning-free robust regression model
\begin{eqnarray}\label{Wilcoxon-score-loss}
\min_{\beta\in\mathcal{R}^{p}}\left\{\frac{1}{n(n-1)}\sum_{1\leq i<j\leq n}\left|(b_{i}-X_{i}\beta)-(b_{j}-X_{j}\beta)\right|+\lambda\|\beta\|_{1}\right\},
\end{eqnarray}
where $X_{i}$ is the $i$-th row of the matrix $X$ and $b_{i}$ is the $i$-th component of the vector $b$.
The loss in (\ref{Wilcoxon-score-loss}) was originally from the classical nonparametric statistics (see e.g., \cite{Hettmansperger1998}) and is equivalent to Jaeckel's dispersion function (see \cite{Jaeckel1972}) with Wilcoxon scores. For simplicity, we call the problem (\ref{Wilcoxon-score-loss}) the tuning-free robust Lasso problem. It has been shown in \cite{WangPeng} that the model (\ref{Wilcoxon-score-loss}) is very close to the Lasso for normal random errors and is robust with substantial efficiency gain for heavy-tailed errors.

As we know, it is computationally challenge to solve the problems of Lasso, square-root Lasso and tuning-free robust Lasso, though they are all convex. Many numerical algorithms have been proposed to solve the Lasso problem, including the accelerated proximal gradient (APG) method \cite{BeckT2009}, the interior-point method (IPM) \cite{KimKLBG2007}, the least angle regression (LARS) \cite{EfronHJT} and the semismooth Newton augmented Lagrangian method (SSNAL) \cite{LiSunToh2018}, and so on. Due to the nonsmooth loss in the objective function, the authors in \cite{Liuhan2015} applied the alternating direction method of multipliers (ADMM) to solve the square-root Lasso problem. In a recent work \cite{TangWangSunToh2020}, an efficient proximal majorization-minimization algorithm (PMM) has been proposed to solve the square-root regression problems with nonconvex regularizers. As for the tuning-free robust  Lasso problem, the authors in \cite{WangPeng} reformulated this problem to an $(n(n-1)+2p)$-dimensional linear programming problem with $(n(n-1)+2p)$ constraints. The well known IPM, which is implemented in some commercial solvers such as Gurobi and Mosek, is a powerful algorithm for solving large-scale problems. When the feature dimension $p$ and/or the sample size $n$ is very large, the computational cost of linear programming is unacceptable.

Although the $\ell_{1}$ regularizer has many attractive properties, for large regression coefficients the shrinkage introduced by the Lasso tends to be significantly biased toward 0 and leads to biased estimate. Many alternative regularizers such as the smoothly clipped absolute deviation (SCAD) regularizer \cite{Fan2001,Fan2004} and the minimax concave penalty (MCP) \cite{Zhang2010} were designed to diminish this bias. All these regularization regression models have the so-called oracle property for individual variables, which means that the corresponding regularized estimator is equal to the least squares estimator assuming the model is known with high probability under appropriate conditions. By making a local linear approximation to the regularizers, Zou and Li \cite{Zouli2008} applied the LARS algorithm to solve these optimization problems. It is known that many of these nonconvex regularizers, which are surrogate sparsity functions, can be expressed as the difference of two convex functions (DC) \cite{Ahn2017,Le-thi2015}. By the DC property of the objective function, we can naturally apply a majorization-minimization algorithm to solve the nonconvex problem.

In this paper, our main contribution is in proposing an efficient algorithm for the tuning-free robust regression problems with nonconvex regularizations, which takes the following form
\begin{eqnarray}\label{tuning-free robust-problem}
\min_{\beta\in\mathcal{R}^{p}}\left\{g(\beta):=\frac{1}{n(n-1)}\sum_{1\leq i<j\leq n}\left|(b_{i}-X_{i}\beta)-(b_{j}-X_{j}\beta)\right|+\lambda q_{1}(\beta)-q_{2}(\beta)\right\},
\end{eqnarray}
where $q_{1}:\mathcal{R}^{p}\rightarrow \mathcal{R}$ is a convex function with a strongly semismooth proximal mapping and
$q_{2}:\mathcal{R}^{p}\rightarrow\mathcal{R}$ is a convex smooth
function (the dependence of $q_{2}$ on $\lambda$ has been dropped here). Since the proximal mappings of many commonly used functions such as the $\ell_{1}$ and $\ell_{2}$ functions are strongly semismooth \cite{Meng2005}, the assumption on $q_{1}$ is mild. For the DC structure of the regularization function, we design a two stage proximal-proximal majorization-minimization (PPMM) algorithm to solve the optimization problem (\ref{tuning-free robust-problem}). For each stage of the PPMM algorithm, the corresponding subproblem is convex and contains two nonsmooth terms, i.e., the loss function and the regularization function $q_{1}$. We apply a semismooth Newton (SSN) method based proximal point algorithm (PPA) to solve the dual problem of each subproblem. Note that the Clarke generalized Jacobian of the subproblem is singular, and using the SSN method directly is not feasible. To overcome this difficulty, we have to apply the SSN method to a sequentially regularized subproblems.

The remaining parts of this paper are organized as follows. In Section \ref{sec:Preliminaries}, we introduce
some basic knowledge which will be used in this paper. In Section \ref{sec:PPA}, we briefly introduce the preconditioned PPA. In Section \ref{sec:PPMM}, we present the details of the PPMM algorithm, and we also prove the convergence and the convergence rate of the proposed algorithm based on the Kurdyka-{\L}ojasiewicz (KL) property. Furthermore, we prove that the algorithm converges to a d-stationary point. In Section \ref{sec:Numerical experiments}, we compare our algorithm with the existing algorithms to demonstrate that the proposed algorithm can solve the tuning-free robust regression problem efficiently.

\subsection{Additional notations}
\label{sec:Additional notations}

Let $\langle\cdot,\cdot\rangle$ and $\|\cdot\|$ be the standard Euclidean inner product
and norm in the space $\mathcal{R}^{n}$. For any $x\in\mathcal{R}^{n}$ and a given self-adjoint positive semidefinite matrix $M\in\mathcal{R}^{n\times n}$, define $\|x\|_{M}:=\sqrt{\langle x,Mx\rangle}$. The largest eigenvalue of $M$ is denoted by $\lambda_{\max}(M)$. Given a set $C\subseteq\mathcal{R}^{n}$, the weighted distance of $x$ to $C$ is defined by $\mbox{dist}_{M}(x,C):=\inf_{y\in C}\{\|y-x\|_{M}\}$. If $C=\emptyset$, we have that $\mbox{dist}_{M}(x,C)=+\infty$ for all $x\in\mathcal{R}^{n}$. Denote the identity matrix of order $n$ by $I_{n}$. If $M$ is an identity matrix, we just omit the subscript matrix $M$. We define the indicator function $\delta_{C}$ of the set $C$ by $\delta_{C}(x)=0$ if $x\in C$, otherwise $\delta_{C}(x)=+\infty$. The intersection of all the convex sets containing $C$ is called the convex hull of $C$ and is denoted by conv$(C)$.

\section{Preliminaries}
\label{sec:Preliminaries}
In this section, we present some basic preliminaries that will be used in this paper.

For any real valued function $r:\mathcal{R}^{n}\rightarrow\overline{\mathcal{R}}$, the conjugate to $r$ is defined by
\begin{eqnarray*}
r^{*}(x):=\sup_{v\in \textrm{dom}(r)}\Big\{\langle x,v\rangle-r(v)\Big\}.
\end{eqnarray*}
We call $r$ a proper function if there exists at least one $x\in\mathcal{R}^{n}$ such that $r(x)<+\infty$ and $r(x)>-\infty$ for all $x\in\mathcal{R}^{n}$, or in other words, if $\mbox{dom}(r)$ is a nonempty set and $r$ is finite; otherwise, it is improper.
For a proper lower semicontinuous function $r$, the Moreau envelope function $e_{\sigma r}$ and the proximal mapping $\mbox{Prox}_{\sigma r}$ with parameter $\sigma>0$ of $r$ are defined respectively by
\begin{eqnarray*}
e_{\sigma r}(x)&:=&\inf_{w\in\mathcal{R}^{n}}\Big\{r(w)+\frac{1}{2\sigma}\|w-x\|^{2}\Big\},\\
\mbox{Prox}_{\sigma r}(x)&:=&\mathop{\mbox{argmin}}_{w\in\mathcal{R}^{n}}\Big\{r(w)+\frac{1}{2\sigma}\|w-x\|^{2}\Big\}.
\end{eqnarray*}
When $r$ is also convex, it is known from Theorem 2.26 of \cite{Rockafellar98} that the Moreau envelope function $e_{\sigma r}$ is convex and continuously differentiable with
\begin{eqnarray*}
\nabla e_{\sigma r}(x)=\frac{1}{\sigma}\big(x-\mbox{Prox}_{\sigma r}(x)\big)
\end{eqnarray*}
and the proximal mapping $\mbox{Prox}_{\sigma r}$ is single-valued and continuous with the following Moreau's identity (see, e.g., Theorem 31.5 of \cite{Rockafellar70}) holds
\begin{eqnarray*}
\mbox{Prox}_{\sigma r}(x)+\sigma\mbox{Prox}_{r^{*}/\sigma}(x/\sigma)=x,\quad \forall\ x\in\mathcal{R}^{n}.
\end{eqnarray*}

Consider a proper function $r:\mathcal{R}^{n}\rightarrow(-\infty,+\infty]$ with $x\in$ $\mbox{dom}(r)$, the regular subdifferential of $r$ at $x$ is defined by
\begin{eqnarray*}
\widehat{\partial}r(x):=\Big\{v\in\mathcal{R}^{n}\ \Big|\ r(x')\geq r(x)+\langle v,x'-x\rangle + o(\|x'-x\|)\Big\},
\end{eqnarray*}
which can be written equivalently as
\begin{eqnarray*}
\widehat{\partial}r(x):=\Big\{v\in\mathcal{R}^{n}\ \Big|\ \mathop{\liminf}_{x\neq x'\rightarrow x}\frac{r(x')-r(x)-\langle v,x'-x\rangle}{\|x'-x\|}\geq 0\Big\}.
\end{eqnarray*}
If $x\notin\mbox{dom}(r)$, $\widehat{\partial}r(x)=\emptyset$.
The limiting subdifferential of $r$ at $x$ is defined by
\begin{eqnarray*}
\partial r(x):=\Big\{v\in\mathcal{R}^{n}\ \Big|\ \exists\ x^{\upsilon}\rightarrow x \ \mbox{and}\ v^{\upsilon}\in\widehat{\partial}r(x^{\upsilon})\ \mbox{with}\ v^{\upsilon}\rightarrow v\Big\}.
\end{eqnarray*}
It is known from Theorem 8.6 of \cite{Rockafellar98} that for $x\in\mbox{dom}(r)$, $\widehat{\partial}r(x)$ and $\partial r(x)$ are closed with $\widehat{\partial}r(x)$ convex and $\widehat{\partial}r(x)\subseteq \partial r(x)$. If $r$ is a proper convex function with $x\in\mbox{dom}(r)$, the regular subdifferential and the limiting subdifferential coincide with the set of subgradients of $r$ at $x$ in the sense of convex analysis.

A function $R:\mathcal{R}^{n}\rightarrow\mathcal{R}^{m}$ is called directionally differentiable at a point $x\in\mathcal{R}^{n}$ in a direction $v\in\mathcal{R}^{n}$ if the limit
\begin{eqnarray*}
R'(x,v):=\lim_{t\downarrow 0}\frac{R(x+tv)-R(x)}{t}
\end{eqnarray*}
exists. If $R$ is directionally differentiable at $x$ in every direction $v\in\mathcal{R}^{n}$, we say that $R$ is directionally differentiable at $x$. If $m=1$, $R'(x;v)$ coincides with the classic directional derivative of a real valued function at $x$ in a direction $v$.

A point $x\in\mathcal{R}^{n}$ is called a directional-stationary (d-stationary) point of a proper function $r$ if $r'(x;v)\geq 0$ holds for any $v\in\mathcal{R}^{n}$.
According to generalized Fermat's rule (see e.g., Theorem 10.1 of \cite{Rockafellar98}), the condition $0\in\widehat{\partial}r(x)$ is a necessary condition for $x$ to be a local minimizer of $r$. Furthermore, if $r$ is directionally differentiable at $x$ and locally Lipschitz continuous near $x$, $0\in\widehat{\partial}r(x)$ if and only if $x$ is a d-stationary point of $r$. If the function $r$ is proper and convex, the condition $0\in\widehat{\partial}r(x)$ is not just necessary for a local minimum but also sufficient for a global minimum.

Let $\mathcal{U}\subseteq\mathcal{R}^{n}$ be open and $R:\mathcal{U}\rightarrow \mathcal{R}^{m}$ be a given vector-valued function that is locally Lipschitz continuous. Let $\mathcal{U}_{R}$ be the subset of $\mathcal{U}$ consisting of the points where $R$ is F(r\'{e}chet)-differentiable and $JR(x)\in\mathcal{R}^{m\times n}$ be the Jacobian matrix of $R$ at $x\in\mathcal{U}_{R}$. It is known from Theorem 9.60 of \cite{Rockafellar98} that the set $\mathcal{U}\setminus\mathcal{U}_{R}$ is negligible. For any $x\in\mathcal{U}$, define the B-subdifferential of $R$ at $x$ by
\begin{eqnarray*}
\partial_{B}R(x):=\Big\{V\in\mathcal{R}^{m\times n}\ \Big|\ \exists\ x^{\upsilon}\rightarrow x\ \mbox{with}\ x^{\upsilon}\in\mathcal{U}_{R}\ \mbox{and}\ JG(x^{\upsilon})\rightarrow V\Big\}.
\end{eqnarray*}
The Clarke subdifferential of $R$ at $x$ is defined by $\partial R(x):=\mbox{conv}(\partial_{B}R(x))$.

For further discussion, we introduce the concept of semismoothness below. See \cite{Mifflin1977,QiSun1993} for more details.
\begin{definition}
Given a locally Lipschitz continuous function $R:\mathcal{U}\subseteq\mathcal{R}^{n}\rightarrow\mathcal{R}^{m}$ with an open set $\mathcal{U}$ and a nonempty, compact valued and upper-semicontinuous multifunction $\mathcal{K}:\mathcal{U}\rightrightarrows\mathcal{R}^{m\times n}$, we call $R$ a semismooth function at $x\in\mathcal{U}$ with respect to the multifunction $\mathcal{K}$ if
\begin{enumerate}[(i)]
\item $R$ is directionally differentiable at $x$;
\item $R(x+\Delta x)-R(x)-\Gamma\Delta x=o(\|\Delta x\|)$, for any $\Gamma\in\mathcal{K}(x+\Delta x)$, $\Delta x\in\mathcal{R}^{n}$ and $\Delta x\rightarrow 0$.
\end{enumerate}
Furthermore, if (ii) is replaced by
\begin{eqnarray*}
R(x+\Delta x)-R(x)-\Gamma\Delta x=O(\|\Delta x\|^{1+\gamma}),\ \mbox{for any}\ \Gamma\in\mathcal{K}(x+\Delta x),\ \Delta x\in\mathcal{R}^{n}\ \mbox{and}\ \Delta x\rightarrow 0,
\end{eqnarray*}
where $\gamma$ is a constant, then we call $R$ a $\gamma$-order (strongly if $\gamma=1$) semismooth function at $x$ with respect to $\mathcal{K}$. The function $R$ is said to be a semismooth function on $\mathcal{U}$ with respect to $\mathcal{K}$ if it is semismooth everywhere on $\mathcal{U}$ with respect to $\mathcal{K}$.
\end{definition}

The Kurdyka-{\L}ojasiewicz (KL) property (see e.g., \cite{bolte-pauwels2016,bolte-sabach-teboulle2014}) plays a central role in our further analysis.
\begin{definition}
A proper lower semicontinuous function $r:\mathcal{R}^{n}\rightarrow(-\infty,+\infty]$ is said to have the KL property at $x\in\mbox{dom}(\partial r)$ if there exist $\eta\in(0,+\infty]$, a neighbour $\mathcal{U}$ of $x$ and a continuous concave function $\varphi:[0,\eta)\rightarrow [0,+\infty)$ satisfying
\begin{enumerate}[$\quad\ $(1)]
\item $\varphi(0)=0$;
\item $\varphi$ is continuous at 0 and continuously differentiable on $(0,\eta)$;
\item$\varphi'(s)>0$, for all $0<s<\eta$
\end{enumerate}
such that the KL inequality $\varphi'(r(x')-r(x))\mbox{dist}(0,\partial r(x'))\geq 1$ holds for any $x'\in \mathcal{U}$ and $r(x)<r(x')<r(x)+\eta$.
If $r$ satisfies the KL property at each point of $\mbox{dom}(\partial r)$, then $r$ is called a KL function.
\end{definition}


\section{The preconditioned PPA}
\label{sec:PPA}
The PPA is a classical optimization algorithm, which dates back to \cite{Martinet1970} and was extensively studied in \cite{Rockafellar1976a,Rockafellar1976b}. The preconditioned PPA which is a generalization of the PPA was first studied by \cite{LiSunToh2020}.
Consider a lower semicontinuous proper convex function $f:\mathcal{R}^{n}\rightarrow (-\infty,+\infty]$. Given a sequence of positive real numbers $\{\sigma_{k}\}$ such that $0<\sigma_{k}\uparrow\sigma_{\infty}\leq +\infty$, $x^{0}\in\mathcal{R}^{n}$ and $\{M_{k}\}$ being a sequence of self-adjoint positive definite matrix in $\mathcal{R}^{n\times n}$ satisfying
\begin{eqnarray*}
(1+\varrho_{k})M_{k}\succeq M_{k+1},\quad M_{k}\succeq\lambda_{\min}I_{n},\ \forall\ k\geq 0, \quad \limsup_{k\rightarrow\infty}\lambda_{\max}(M_{k})=\lambda_{\infty}
\end{eqnarray*}
with $\{\varrho_{k}\}$ a nonnegative summable sequence and $0<\lambda_{\min}\leq\lambda_{\infty}<+\infty$,
the preconditioned PPA generates a sequence $\{x^{k}\}$ approximately by
\begin{eqnarray*}
x^{k+1}\approx\bar{x}^{k+1}=\mathop {\rm argmin }_{x\in\mathcal{R}^{n}}\left\{f(x)+\frac{1}{2\sigma_{k}}\|x-x^{k}\|_{M_{k}}^{2}\right\}.
\end{eqnarray*}
If $M_{k}\equiv I_{n}$ for all $k\geq 0$, the preconditioned PPA becomes to the classical PPA.
There are two general criteria for the approximate computation of $x^{k+1}$,
\begin{eqnarray*}
&&(A)\quad \|x^{k+1}-\bar{x}^{k+1}\|_{M_{k}}\leq \varepsilon_{k},\quad0\leq\varepsilon_{k},\quad \sum_{k=0}^{\infty}\varepsilon_{k}<+\infty,\\
&&(B)\quad \|x^{k+1}-\bar{x}^{k+1}\|_{M_{k}}\leq \delta_{k}\|x^{k+1}-x^{k}\|_{M_{k}},\quad 0\leq\delta_{k}<1,\quad \sum_{k=0}^{\infty}\delta_{k}<+\infty.
\end{eqnarray*}

A multifunction $\mathcal{F}:\mathcal{R}^{d}\rightrightarrows\mathcal{R}^{d}$ is said to be locally upper Lipschitz continuous at $x$ with modulus $\kappa$ which is independent of $x$ if there exists a neighbourhood $\mathcal{U}$ of $x$ such that $\mathcal{F}(y)\subseteq\mathcal{F}(x)+\kappa\|y-x\|\mathbb{B}_{d}$, $\forall\ y\in\mathcal{U}$, where $\mathbb{B}_{d}$ is a unit ball in $\mathcal{R}^{d}$. A multifunction $\mathcal{F}$ is said to be piecewise polyhedral if its graph $\mbox{gph}\mathcal{F}:=\{(x,y)\ |\ y\in\mathcal{F}(x)\}$ is the union of finitely many polyhedral convex sets. The inverse of a piecewise polyhedral multifunction is also piecewise polyhedral. It has been shown by \cite{Robinson1981} that a piecewise polyhedral set-valued mapping has a fundamental locally upper Lipschitz continuous property.
\begin{proposition}
Let $\mathcal{F}:\mathcal{R}^{d}\rightrightarrows\mathcal{R}^{d}$ be a piecewise polyhedral set-valued mapping. For any $x\in\mathcal{R}^{d}$, there exists a constant $\kappa$ independent of $x$ such that $\mathcal{F}$ is locally upper Lipschitz continuous with modulus $\kappa$.
\end{proposition}

We review some convergence results of the preconditioned PPA; one may also refer to \cite{LiSunToh2020}.
\begin{theorem}
\label{thm:convergence of PPA}
Suppose that $\Omega:=\{x\ |\ 0\in\partial f(x)\}\neq\emptyset$. Let $\{x^{k}\}$ be any sequence generated by the preconditioned PPA under the criterion $(A)$. Then the sequence $\{x^{k}\}$ is bounded and
\begin{eqnarray*}
\mbox{dist}_{M_{k+1}}(x^{k+1},\Omega)\leq(1+\varrho_{k})\mbox{dist}_{M_{k}}(x^{k},\Omega)+(1+\varrho_{k})\varepsilon_{k},\quad \forall\ k\geq 0.
\end{eqnarray*}
In addition, the sequence $\{x^{k}\}$ converges to a point $x^{\infty}$ with $0\in\partial f(x^{\infty})$.
\end{theorem}

\begin{assumption}\label{assump-error-bound}
The operator $\partial f$ satisfies the following error bound condition. For any $\delta>0$, there exists $\kappa>0$ such that
\begin{eqnarray*}
\mbox{dist}(x,(\partial f)^{-1}(0))\leq\kappa\mbox{dist}(0,\partial f(x)),\quad \forall\,x\in\{x\,|\,\mbox{dist}(x,(\partial f)^{-1}(0))\leq \delta\}.
\end{eqnarray*}
\end{assumption}

\begin{remark}
The error bound condition is critical to obtain the convergence rate of the preconditioned PPA. It has been proven by Lemma 2.4 of \cite{LiSunToh2020} that a multifunction satisfies the error bound condition if its inverse is locally upper Lipschitz continuous at the origin.
\end{remark}

\begin{theorem}
\label{thm:convergence rate of PPA}
Assume that $\Omega\neq\emptyset$ and Assumption \ref{assump-error-bound} holds.
Let $t$ be a positive number satisfying $t>\sum\limits_{k=0}^{\infty}\varepsilon_{k}(1+\varrho_{k})$ and the initial point $x^{0}$ be a point such that
\begin{eqnarray*}
\mbox{dist}_{M_{0}}(x^{0},\Omega)\leq\frac{t-\sum\limits_{k=0}^{\infty}\varepsilon_{k}(1+\varrho_{k})}{\prod\limits_{k=0}^{\infty}(1+\varrho_{k})}.
\end{eqnarray*}
Let $\{x^{k}\}$ be any sequence generated by the preconditioned PPA under the criteria $(A)$ and $(B)$ with $\{\sigma_{k}\}$ nondecreasing. Then we have
\begin{eqnarray*}
\mbox{dist}_{M_{k+1}}(x^{k+1},\Omega)\leq\theta_{k}\mbox{dist}_{M_{k}}(x^{k},\Omega),\ \mbox{for all}\ k\geq 0,
\end{eqnarray*}
where
\begin{eqnarray*}
\theta_{k}=(1+\varrho_{k})(1-\delta_{k})^{-1}\left(\delta_{k}+\frac{(1+\delta_{k})\kappa\lambda_{max}(M_{k})}{\sqrt{\sigma_{k}^{2}+\kappa^{2}\lambda_{max}^{2}(M_{k})}}\right)
\end{eqnarray*}
and
\begin{eqnarray*}
\limsup_{k\rightarrow\infty}\theta_{k}=\theta_{\infty}=\frac{\kappa\lambda_{\infty}}{\sqrt{\sigma_{\infty}^{2}+\kappa^{2}\lambda_{\infty}^{2}}}<1\quad\mbox{with}\ \theta_{\infty}=0\ \mbox{if}\ \sigma_{\infty}=\infty.
\end{eqnarray*}
\end{theorem}

\section{The PPMM algorithm}
\label{sec:PPMM}
In this section, we introduce the PPMM algorithm to solve the nonconvex tuning-free robust regression problem. Firstly, we obtain an initial point to warm start the PPMM algorithm by solving a convex relaxation problem with omitting the term $-q_{2}$ of the problem (\ref{tuning-free robust-problem}). We solve this subproblem by the SSN method based PPA. Secondly, the solution of the nonconvex problem is obtained by solving a series of majorizied problems. Specifically, the majorized problem is constructed by linearizing the concave term $-q_{2}$ at the current iterate point $\tilde{\beta}$. It is essentially important to solve those subproblems efficiently and accurately. However, it is not easy to accomplish this task. The main difficulty is that the Clarke generalized Jacobian of the subproblem may be probably singular. To overcome this difficulty, we introduce a special proximal term $\frac{\tau}{2}\|\beta-\tilde{\beta}\|^{2}+\frac{\iota}{2}\|X(\beta-\tilde{\beta})\|^{2}$ ($\tau>0$, $\iota>0$) and sequentially solve these regularized problems. After introducing the proximal term, the dual of the problem can be written explicitly as an unconstrained convex problem that can be solved by the SSN method.

For $\tilde{\beta},\tilde{w}\in\mathcal{R}^{p}$, $\tau>0$ and $\iota>0$, denote
\begin{eqnarray}\label{subproblem}
f(\beta;\tau,\iota,\tilde{\beta},\tilde{w})&:=&\frac{1}{n(n-1)}\sum_{1\leq i<j\leq n}\left|(b_{i}-X_{i}\beta)-(b_{j}-X_{j}\beta)\right|+\lambda q_{1}(\beta)\nonumber\\
&&-q_{2}(\tilde{\beta})-\langle\tilde{w},\beta-\tilde{\beta}\rangle+\frac{\tau}{2}\|\beta-\tilde{\beta}\|^{2}+\frac{\iota}{2}\|X(\beta-\tilde{\beta})\|^{2}.
\end{eqnarray}
Now we formally present the PPMM algorithm for the tuning-free robust regression problem.

\bigskip

\ni\fbox{\parbox{\textwidth}{\noindent{\bf Algorithm 1 (PPMM): \label{alg:PPMM}}
\begin{description}
\item [Stage 1:] (\textbf{The preconditioned PPA with} $M_{k}\equiv I_{p}+X^{T}X$) Given $\tau_{1,0}>0$, $\beta^{0,0}\in\mathcal{R}^{p}$, set $k=0$, iterate:
\item [Step 1.1.] Apply Algorithm 2 to find an approximate solution $u^{0,k+1}$ of the dual problem of the following optimization problem
\begin{eqnarray*}
\min_{\beta\in\mathcal{R}^{p}}\Big\{f(\beta;\tau_{1,k},\tau_{1,k},\beta^{0,k},0)\Big\}.
\end{eqnarray*}
Then $\beta^{0,k+1}=\mbox{Prox}_{\tau_{1,k}^{-1}\lambda q_{1}}(-\tau_{1,k}^{-1}X^{T}u^{0,k+1}+\beta^{0,k})$.
If a desired stopping criterion is satisfied, let $\beta^{0}=\beta^{0,k+1}$, go to Stage 2; otherwise go to Step 1.2.
\item [Step 1.2.]  Update $\tau_{1,k+1}=\rho_{1,k}\tau_{1,k}$, $\rho_{1,k}\in(0,1)$. Set $k:=k+1$ and return to Step 1.1.
\item [Stage 2:]  Given $\tau_{2,0}>0$ and $\iota_{2,0}>0$, set $k=0$, iterate:
\item [Step 2.1.] Apply Algorithm 2 to find an approximate solution $u^{k+1}$ of the dual problem of the following optimization problem
\begin{eqnarray}\label{PPMM-subproblem}
\min_{\beta\in\mathcal{R}^{p}}\Big\{f(\beta;\tau_{2,k},\iota_{2,k},\beta^{k},\nabla q_{2}(\beta^{k}))\Big\}.
\end{eqnarray}
Then $\beta^{k+1}=\mbox{Prox}_{\tau_{2,k}^{-1}\lambda q_{1}}(-\tau_{2,k}^{-1}X^{T}u^{k+1}+\tau_{2,k}^{-1}\nabla q_{2}(\beta^{k})+\beta^{k})$.
If a desired stopping criterion is satisfied, terminate; otherwise go to Step 2.2.
\item [Step 2.2.]  Update $\tau_{2,k+1}=\rho_{2,k}\tau_{2,k}$, $\iota_{2,k+1}=\rho'_{2,k}\iota_{2,k}$, $\rho_{2,k}\in(0,1),\rho'_{2,k}\in(0,1)$. Set $k:=k+1$ and return to Step 2.1.
\end{description}}}

\bigskip

\subsection{The subproblem of the PPMM algorithm}
Given $\tau>0$, $\iota>0$ and $\tilde{\beta}\in\mathcal{R}^{p}$, consider the following optimization problem
\begin{eqnarray}\label{subproblem}
\min_{\beta\in\mathcal{R}^{p}}\Big\{f(\beta;\tau,\iota,\tilde{\beta},\tilde{w})\Big\}.
\end{eqnarray}
Denote $h(y) := \frac{1}{n(n-1)}\sum\limits_{1\leq i<j\leq n}\left|y_{i}-y_{j}\right|$, the problem (\ref{subproblem}) is equivalent to the following minimization problem
\begin{eqnarray}\label{subproblem-1}
\min_{{\beta\in\mathcal{R}^{p}}\atop{y\in\mathcal{R}^{n}}}\Big\{h(y)+\lambda q_{1}(\beta)-\langle\tilde{w},\beta-\tilde{\beta}\rangle+\frac{\tau}{2}\|\beta-\tilde{\beta}\|^{2}+\frac{\iota}{2}\|y+b-X\tilde{\beta}\|^{2}\ \Big|\ X\beta-y-b=0\Big\}.
\end{eqnarray}
We can also write the dual problem of (\ref{subproblem-1}) as
\begin{eqnarray}\label{subproblem-dual}
&&\min_{u\in\mathcal{R}^{n}}\Big\{\phi(u):=\frac{\iota}{2}\|\iota^{-1}u+X\tilde{\beta}-b\|^{2}-\frac{\iota}{2}\|\iota^{-1}u+X\tilde{\beta}-b-\mbox{Prox}_{\iota^{-1}h}(\iota^{-1}u+X\tilde{\beta}-b)\|^{2}\nonumber\\
&&\quad\quad\quad-h(\mbox{Prox}_{\iota^{-1}h}(\iota^{-1}u+X\tilde{\beta}-b))+\frac{\tau}{2}\|-\tau^{-1}X^{T}u+\tau^{-1}\tilde{w}+\tilde{\beta}\|^{2}\\
&&\quad\quad\quad-\frac{\tau}{2}\|-\tau^{-1}X^{T}u+\tau^{-1}\tilde{w}+\tilde{\beta}-\mbox{Prox}_{\tau^{-1}\lambda q_{1}}(-\tau^{-1}X^{T}u+\tau^{-1}\tilde{w}+\tilde{\beta})\|^{2}\nonumber\\
&&\quad\quad\quad-\lambda q_{1}(\mbox{Prox}_{\tau^{-1}\lambda q_{1}}(-\tau^{-1}X^{T}u+\tau^{-1}\tilde{w}+\tilde{\beta}))+\langle u,b\rangle\Big\}.\nonumber
\end{eqnarray}
Let $\bar{u}$ be a solution of the problem (\ref{subproblem-dual}), then we have $\bar{\beta}=\mbox{Prox}_{\tau^{-1}\lambda q_{1}}(-\tau^{-1}X^{T}u+\tau^{-1}\tilde{w}+\tilde{\beta})$ and $\bar{y}=\mbox{Prox}_{\iota^{-1}h}(\iota^{-1}u+X\tilde{\beta}-b)$, which is a solution of the primal problem (\ref{subproblem-1}).

The objective function of the problem (\ref{subproblem-dual}) is smooth and finding the optimal solution of the problem (\ref{subproblem-dual}) is equivalent to solving the following system of equations
\begin{eqnarray}\label{subproblem-grad}
\nabla\phi(u)=\mbox{Prox}_{\iota^{-1}h}(\iota^{-1}u+X\tilde{\beta}-b)-X\mbox{Prox}_{\tau^{-1}\lambda q_{1}}(-\tau^{-1}X^{T}u+\tau^{-1}\tilde{w}+\tilde{\beta})+b=0.
\end{eqnarray}
One may apply the SSN method to solve the aforementioned problem (\ref{subproblem-grad}). However, the Clarke generalized Jacobian of the problem (\ref{subproblem-grad}) may be probably singular and the SSN method cannot be directly applied. To overcome this difficulty, we deal with the problem (\ref{subproblem-dual}) in an alternative way by applying the PPA with the inner subproblem solved by the SSN method.

Given $u^{0}\in\mathcal{R}^{n}$, let $\{u^{i}\}$ be the sequence generated by the PPA. We adopt
\begin{eqnarray*}
(B'')\quad\|\nabla\phi(u^{i+1})+\sigma_{i}^{-1}(u^{i+1}-u^{i})\|\leq\frac{\delta_{i}}{\sigma_{i}}\|u^{i+1}-u^{i}\|,\ 0\leq\delta_{i}<1,\ \sum_{i=0}^{\infty}\delta_{i}<\infty
\end{eqnarray*}
as the stopping criterion of each iteration.  We apply the SSN method to compute an approximate solution of each subproblem. Let $\hat{\phi}(u;\sigma_{i},u^{i}):=\phi(u)+\frac{1}{2\sigma_{i}}\|u-u^{i}\|^{2}$. Since the mappings $\mbox{Prox}_{\iota^{-1}h}(\cdot)$ and $\mbox{Prox}_{\tau^{-1}\lambda q_{1}}(\cdot)$ are both Lipschitz continuous, the following multifunction
\begin{eqnarray*}
\hat{\partial}^{2}\hat{\phi}(u;\sigma_{i},u^{i})&:=&\iota^{-1}\partial\mbox{Prox}_{\iota^{-1}h}(\iota^{-1}u+X\tilde{\beta}-b)\\&&+\tau^{-1}X\partial\mbox{Prox}_{\tau^{-1}\lambda q_{1}}(-\tau^{-1}X^{T}u+\tau^{-1}\tilde{w}+\tilde{\beta})X^{T}
+\sigma_{i}^{-1}I_{n}
\end{eqnarray*}
is well defined. It is known from \cite{HiriartSN1984} that
\begin{eqnarray*}
\partial^{2}\hat{\phi}(u;\sigma_{i},u^{i})(d)=\hat{\partial}^{2}\hat{\phi}(u;\sigma_{i},u^{i})(d),\quad \forall\ d\in\mathcal{R}^{n},
\end{eqnarray*}
where $\partial^{2}\hat{\phi}(u;\sigma_{i},u^{i})$ is the Clarke generalized Jacobian of $\hat{\phi}(u;\sigma_{i},u^{i})$ at $u$.
Let us select $U\in\partial\mbox{Prox}_{\iota^{-1}h}(\iota^{-1}u+X\tilde{\beta}-b)$ and $V\in\partial\mbox{Prox}_{\tau^{-1}\lambda q_{1}}(-\tau^{-1}X^{T}u+\tau^{-1}\tilde{w}+\tilde{\beta})$, we obtain $H:=\iota^{-1}U+\tau^{-1}XVX^{T}+\sigma_{i}^{-1}I_{n}\in\partial^{2}\hat{\phi}(u;\sigma_{i},u^{i})$.

Now we introduce the PPA and the SSN method as below.
\bigskip

\ni\fbox{\parbox{\textwidth}{\noindent{\bf Algorithm 2 (PPA): \label{alg:PPA}} Given $\sigma_{0}>0$, choose $u^{0}\in\mathcal{R}^{n}$. Set $i=0$ and iterate:
\begin{description}
\item [Step 1.] Apply Algorithm 3 to solve the following problem
\begin{eqnarray}
u^{i+1}\approx\mathop{\rm argmin}_{u\in\mathcal{R}^{n}}\hat{\phi}(u;\sigma_{i},u^{i}).\label{eq:PPA}
\end{eqnarray}
\item [Step 2.] If a desired stopping criterion is satisfied, terminate; otherwise, update $\sigma_{i+1}=\rho_{i}\sigma_{i}$ with $\rho_{i}\geq1$ and return to Step 1.
\end{description}}}

\bigskip

\bigskip

\ni\fbox{\parbox{\textwidth}{\noindent{\bf Algorithm 3 (SSN): \label{alg:Newton-CG1}} Given $\iota,\tau,\sigma>0,$ $\tilde{\beta},\tilde{w}\in\mathcal{R}^{p}$,
$\mu\in(0,\frac{1}{2}), \ \overline{\eta}\in(0,1),
\varsigma\in(0,1],\textrm{
  and }
  \delta\in(0,1)$, choose $u^{0}\in\mathcal{R}^{n}$. Set $j=0$ and iterate:
\begin{description}
\item [Step 1.] Choose $U^{j}\in\partial\mbox{Prox}_{\iota^{-1}h}(\iota^{-1}u^{j}+X\tilde{\beta}-b)$ and $V^{j}\in\partial\mbox{Prox}_{\tau^{-1}\lambda q_{1}}(-\tau^{-1}X^{T}u^{j}+\tau^{-1}\tilde{w}+\tilde{\beta})$.  Let
$H^{j}=\iota^{-1}U^{j}+\tau^{-1}XV^{j}X^{T}+\sigma^{-1}I_{n}$ and find the
solution $u^{j}$ of the following linear system
\begin{eqnarray*}
H^{j}\Delta u=-\nabla\hat{\phi}(u^{j};\sigma,u^{0})
\end{eqnarray*}
\hspace {22mm} such that
\begin{eqnarray}
\|H^{j}\Delta u^{j}+\nabla\hat{\phi}(u^{j};\sigma,u^{0})\|\leq
\eta_{j}:=\min(\overline{\eta},\|\nabla\hat{\phi}(u^{j};\sigma,u^{0})\|^{1+\varsigma}).\label{ineq:stopping
criterion 1}
\end{eqnarray}
\item [Step 2.]  Set $\alpha_{j}=\delta^{m_{j}}$, where $m_{j}$ is the first nonnegative integer $m$ for which
\begin{eqnarray*}
&\hat{\phi}(u^{j}+\delta^{m}\Delta u^{j};\sigma,u^{0})\leq
\hat{\phi}(u^{j};\sigma,u^{0})+\mu\delta^{m}\langle\nabla
 \hat{\phi}(u^{j};\sigma,u^{0}),(\Delta u^{j})\rangle.&
\end{eqnarray*}
\item [Step 3.]  Set $u^{j+1}=u^{j}+\alpha_{j}\Delta u^{j}$.
\end{description}}}

\bigskip

Since $\mbox{Prox}_{\iota^{-1}h}(\cdot)$ is a piecewise affine function, it is strongly semismooth due to Proposition 7.4.7 of \cite{Facchinei2003}. We have assumed that $\mbox{Prox}_{q_{1}}(\cdot)$ is strongly semismooth. Therefore the function $\nabla\hat{\phi}(\cdot;\sigma,u^{0})$ is strongly semismooth. Hence we can establish the result about the local convergence rate of the SSN method. We just present the main result here without proof.
\begin{theorem}
Suppose that $\sigma<+\infty$. The sequence $\{u^{j}\}$ generated by the SSN method converges to the unique solution $\bar{u}$ of $\nabla\hat{\phi}(v;\sigma,u^{0})=0$ with
\begin{eqnarray*}
\|u^{j+1}-\bar{u}\|=\mathcal{O}(\|u^{j}-\bar{u}\|^{1+\varsigma}).
\end{eqnarray*}
\end{theorem}

\subsection{The computation of $\mbox{Prox}_{\iota^{-1}h}(x)$ and its generalized Jacobian}
In the implementation of the SSN method, it is essentially important to compute the structure of $\partial\mbox{Prox}_{\iota^{-1}h}(x)$. In this subsection, we introduce some basic results about the proximal mapping of $h$ and its corresponding generalized Jacobian which have been considered by \cite{Linmeixia2019}. Let $\rho:=\frac{1}{\iota n(n-1)}$, the proximal mapping related to $\iota^{-1}h$ can be computed by the following proposition.
\begin{proposition}
Let $x\in\mathcal{R}^{n}$ be an arbitrarily given vector and $P_{x}\in\mathcal{R}^{n\times n}$ be a corresponding permutation matrix such that $x^{\downarrow}=P_{x}x$ and $x^{\downarrow}_{1}\geq x^{\downarrow}_{2}\geq\cdots\geq x^{\downarrow}_{n}$. Then the proximal mapping related to $\iota^{-1}h$ can be computed by
\begin{eqnarray*}
\mbox{Prox}_{\iota^{-1}h}(x)=P_{x}^{T}\Pi_{D}(x^{\downarrow}-\rho w),
\end{eqnarray*}
where $w\in\mathcal{R}^{n}$ with $w_{i}=n-2i+1$, $i=1,2,\ldots,n$, $D=\{x\in\mathcal{R}^{n}\ |\ Bx\geq 0\}$ with $B$ being a matrix such that $Bx=[x_{1}-x_{2};\ldots;x_{n-1}-x_{n}]\in\mathcal{R}^{n-1}$ and $\Pi_{D}(\cdot)$ is the metric projection onto $D$ which can be computed by the pool-adjacent-violators algorithm (see \cite{Best1990} for more details).
\end{proposition}

Now we introduce the generalized Jacobian of $\mbox{Prox}_{\iota^{-1}h}(\cdot)$, denoted by $\mathcal{M}(\cdot)$. It is known that $\mbox{Prox}_{\iota^{-1}h}(\cdot)$ is strongly semismooth with respect to $\mathcal{M}(\cdot)$. The details of $\mathcal{M}(\cdot)$ can be found in \cite{Linmeixia2019}.

For a given vector $y\in\mathcal{R}^{n}$, the active index set is denoted by
\begin{eqnarray*}
\mathcal{I}_{A}(y):=\Big\{i\ \Big|\ B_{i}\Pi_{D}(y)=0,\ i=1,2,\ldots,n-1\Big\},
\end{eqnarray*}
where $B_{i}$ is the $i$-th row of the matrix $B$. Based on the active index set $\mathcal{I}_{A}(x^{\downarrow}-\rho w)$, we define a diagonal matrix $\Sigma\in\mathcal{R}^{(n-1)\times(n-1)}$ with each diagonal element taking the following value
\begin{eqnarray*}
\Sigma_{ii}=\left\{\begin{array}{cc}1,& i\in \mathcal{I}_{A}(x^{\downarrow}-\rho w),\\0,& \mbox{otherwise,}\end{array}\right.\quad i=1,2,\ldots,n-1.
\end{eqnarray*}

Now we are ready to get the following proposition which proposes a way to find a computational element of $\mathcal{M}(\cdot)$.
\begin{proposition}
For a given vector $x\in\mathcal{R}^{n}$, we have
\begin{eqnarray*}
P_{x}^{T}\big(I_{n}-B^{T}(\Sigma BB^{T}\Sigma)^{\dagger}B\big)P_{x}\in\mathcal{M}(x),
\end{eqnarray*}
where $(\cdot)^{\dagger}$ denotes the pseduoinverse.
\end{proposition}
\subsection{The convergence analysis of the PPMM algorithm}
For Stage 1 of the PPMM algorithm, it solves a convex and polyhedral problem, hence the convergence results follow from Theorem \ref{thm:convergence of PPA} and \ref{thm:convergence rate of PPA}.

Now we focus on the convergence analysis of Stage 2 of the PPMM algorithm on the nonconvex problem. Let $r^{k+1}$ be the residual vector of the gradient for the objective function of the dual problem related to the problem (\ref{PPMM-subproblem}) at $u^{k+1}$. By the stopping criterion ($B''$) of the PPA, we know that $r^{k+1}$ converges to zero as $\sigma_{i}\rightarrow\infty$. It follows that $r^{k+1}$ satisfies
\begin{eqnarray}\label{cond:algm-ppmm}
2h(r^{k+1})+\iota_{2,k}\|r^{k+1}\|^{2}\leq\frac{\iota_{2,k}}{2}\|X(\beta^{k+1}-\beta^{k})\|^{2},
\end{eqnarray}
if $\sigma_{i}$ is large enough.
For simplicity, denote
\begin{eqnarray*}
f_{k}(\beta)&:=&h(X\beta-b+r^{k+1})+\lambda q_{1}(\beta)-q_{2}(\beta^{k})-\langle\nabla q_{2}(\beta^{k}),\beta-\beta^{k}\rangle\\
&&+\frac{\tau_{2,k}}{2}\|\beta-\beta^{k}\|^{2}+\frac{\iota_{2,k}}{2}\|X(\beta-\beta^{k})+r^{k+1}\|^{2}.
\end{eqnarray*}
Then we can rewrite Stage 2 of the PPMM algorithm as
$$
\beta^{k+1}=\mathop{\mbox{argmin}}\limits_{\beta\in\mathcal{R}^{p}}\Big\{f_{k}(\beta)\Big\}.
$$
 We first introduce the following descent property of the series $\{g(\beta^{k})\}$.
\begin{lemma}\label{lemma-g-descent}
Let $\{\beta^{k}\}$ be a sequence generated by the PPMM algorithm. We have the following descent property
\begin{eqnarray*}
g(\beta^{k})\geq g(\beta^{k+1})+\frac{\tau_{2,k}}{2}\|\beta^{k+1}-\beta^{k}\|^{2}.
\end{eqnarray*}
\end{lemma}
\begin{proof}
Due to the following property of the function $h$
\begin{eqnarray}\label{Lipschitz-h}
|h(x)-h(y)|\leq h(x-y),\ \mbox{for}\ \forall\ x,y\in\mathcal{R}^{n},
\end{eqnarray}
we have
\begin{eqnarray*}
f_{k}(\beta^{k})\leq g(\beta^{k})+h(r^{k+1})+\frac{\iota_{2,k}}{2}\|r^{k+1}\|^{2}.
\end{eqnarray*}
On the other hand, since $q_{2}$ is a convex function and $\beta^{k+1}$ is a minimizer of the function $f_{k}$, we obtain
\begin{eqnarray*}
f_{k}(\beta^{k})&\geq& f_{k}(\beta^{k+1})\\
&\geq& h(X\beta^{k+1}-b+r^{k+1})+\lambda q_{1}(\beta^{k+1})-q_{2}(\beta^{k+1})+\frac{\tau_{2,k}}{2}\|\beta^{k+1}-\beta^{k}\|^{2}\\
&&+\frac{\iota_{2,k}}{2}\|X(\beta^{k+1}-\beta^{k})+r^{k+1}\|^{2}\\
&\geq& g(\beta^{k+1})+\frac{\tau_{2,k}}{2}\|\beta^{k+1}-\beta^{k}\|^{2}+\frac{\iota_{2,k}}{2}\|X(\beta^{k+1}-\beta^{k})\|^{2}-h(r^{k+1})-\frac{\iota_{2,k}}{2}\|r^{k+1}\|^{2}.
\end{eqnarray*}
Combining those above two results, we derive
\begin{eqnarray*}
g(\beta^{k})\geq g(\beta^{k+1})+\frac{\tau_{2,k}}{2}\|\beta^{k+1}-\beta^{k}\|^{2}+\frac{\iota_{2,k}}{2}\|X(\beta^{k+1}-\beta^{k})\|^{2}-2h(r^{k+1})-\iota_{2,k}\|r^{k+1}\|^{2}.
\end{eqnarray*}
By the condition (\ref{cond:algm-ppmm}), the desired result follows.
\end{proof}
The following lemma is similar to Lemmas 5 and 6 in \cite{Cui2018}, and we present it here without proof.
\begin{lemma}\label{lemma-equivalence}
Let $\bar{\beta}\in\mathcal{R}^{p}$. Then $\bar{\beta}$ is a d-stationary point of (\ref{tuning-free robust-problem}) if and only if there exist $\tau,\iota\geq 0$ such that $\bar{\beta}$ is a solution of the following minimization problem
\begin{eqnarray*}
\min_{\beta\in\mathcal{R}^{p}}\Big\{f(\beta;\tau,\iota,\bar{\beta},\nabla q_{2}(\bar{\beta}))\Big\}.
\end{eqnarray*}
\end{lemma}

Now we introduce the convergence result of the PPMM algorithm.
\begin{theorem}\label{thm-convergence}
Assume that the objective function in (\ref{tuning-free robust-problem}) is bounded below and $\{\tau_{2,k}\}$ and $\{\iota_{2,k}\}$ are positive convergent sequences. Then every cluster point of the sequence $\{\beta^{k}\}$ generated by the PPMM algorithm is a d-stationary point of (\ref{tuning-free robust-problem}).
\end{theorem}
\begin{proof}
It is known from Lemma \ref{lemma-g-descent} that $\{g(\beta^{k})\}$ is a non-increasing sequence. Therefore, the sequences $\{g(\beta^{k})\}$ and $\{\|\beta^{k+1}-\beta^{k}\|\}$ both converge and $\lim\limits_{k\rightarrow 0}\|\beta^{k+1}-\beta^{k}\|=0$ due to the assumption of the objective function. Let $\{\beta^{k}\}_{k\in\mathcal{K}}$ be a subsequence of $\{\beta^{k}\}$ with $\lim\limits_{k(\in\mathcal{K})\rightarrow\infty}\beta^{k}=\beta^{\infty}$. We will prove that $\beta^{\infty}$ is a d-stationary point of (\ref{tuning-free robust-problem}).

By the definition of $\beta^{k+1}$, we have
\begin{eqnarray*}
&&f(\beta;\tau_{2,k},\iota_{2,k},\beta^{k},\nabla q_{2}(\beta^{k}))+h(r^{k+1})+\frac{\iota_{2,k}}{2}\|r^{k+1}\|^{2}\\
&\geq& f_{k}(\beta)\geq f_{k}(\beta^{k+1})\\
&\geq& f(\beta^{k+1};\tau_{2,k},\iota_{2,k},\beta^{k},\nabla q_{2}(\beta^{k}))-h(r^{k+1})-\frac{\iota_{2,k}}{2}\|r^{k+1}\|^{2},\quad \forall\ \beta\in\mathcal{R}^{p}.
\end{eqnarray*}
Let $\tau_{\infty}=\lim\limits_{k\rightarrow\infty}\tau_{2,k}$ and $\iota_{\infty}=\lim\limits_{k\rightarrow\infty}\iota_{2,k}$. Since $\|r^{k}\|\rightarrow 0$ due to (\ref{cond:algm-ppmm}) and $h$ is continuous, by letting $k(\in\mathcal{K})\rightarrow\infty$ we derive
\begin{eqnarray*}
f(\beta;\tau_{\infty},\iota_{\infty},\beta^{\infty},\nabla q_{2}(\beta^{\infty}))\geq f(\beta^{\infty};\tau_{\infty},\iota_{\infty},\beta^{\infty},\nabla q_{2}(\beta^{\infty})).
\end{eqnarray*}
Therefore, $\beta^{\infty}\in\mathop{\mbox{argmin}}\limits_{\beta\in\mathcal{R}^{p}}\Big\{f(\beta;\tau_{\infty},\iota_{\infty},\beta^{\infty},\nabla q_{2}(\beta^{\infty}))\Big\}$ and the desired result follows from Lemma \ref{lemma-equivalence}.
\end{proof}

Based on an isolation assumption of the cluster point or the KL property assumption, we establish the local convergence rate of the sequence $\{\beta^{k}\}$ in the following theorem.
\begin{theorem}
Assume that the objective function in (\ref{tuning-free robust-problem}) is bounded below. Let $\{\beta^{k}\}$ be the sequence generated by the PPMM algorithm with $\mathcal{C}^{\infty}$ as the set of all its cluster points. The whole sequence $\{\beta^{k}\}$ will converge to an element of $\mathcal{C}^{\infty}$, if one of the following conditions holds:
\begin{enumerate}[(I)]
\item The set $\mathcal{C}^{\infty}$ contains an isolated element;
\item $\{\beta^{k}\}$ is a bounded sequence. The objective function $g$ in (\ref{tuning-free robust-problem}) has the KL property at $\bar{\beta}$ and $\nabla q_{2}$ is locally Lipschitz continuous near $\bar{\beta}$ for any $\bar{\beta}\in\mathcal{C}^{\infty}$.
\end{enumerate}
Furthermore, based on the condition (2), let $\lim\limits_{k\rightarrow\infty}\beta^{k}=\beta^{\infty}\in\mathcal{C}^{\infty}$ and the function $g$ satisfies the KL property at $\beta^{\infty}$ with an exponent $\gamma\in[0,1)$, we have
\begin{enumerate}[(a)]
\item The sequence $\{\beta^{k}\}$ converges in a finite number of steps, if $\gamma = 0$;
\item The sequence $\{\beta^{k}\}$ converges R-linearly, i.e., there exist $\mu>0$ and $\theta\in[0,1)$ such that $\|\beta^{k}-\beta^{\infty}\|\leq\mu\theta^{^{k}}$, if $0<\gamma\leq\frac{1}{2}$ and $k$ is sufficiently large;
\item The sequence $\{\beta^{k}\}$ converges R-sublinearly, i.e., there exists $\mu>0$ such that $\|\beta^{k}-\beta^{\infty}\|\leq\mu k^{-\frac{1-\gamma}{2\gamma-1}}$, if $\frac{1}{2}<\gamma<1$ and $k$ is sufficiently large.
\end{enumerate}
\end{theorem}
\begin{proof}
It is already known from the proof of Theorem \ref{thm-convergence} that $\lim\limits_{k\rightarrow\infty}\|\beta^{k+1}-\beta^{k}\|=0$. Under the condition (I), the sequence converges to the isolated element of $\mathcal{C}^{\infty}$ due to Proposition 8.3.10 of \cite{Facchinei2003}. In order to prove the convergence rate of the sequence $\{\beta^{k}\}$ based on the condition (II), we first describe some properties as below.
\begin{enumerate}[(i)]
\item The sequence $\{g(\beta^{k})\}$ is descent, i.e., $g(\beta^{k+1})\geq g(\beta^{k})+\frac{\tau_{2,k}}{2}\|\beta^{k+1}-\beta^{k}\|^{2}$;
\item There exists a subsequence $\{\beta^{k}\}_{k\in\mathcal{K}}$ of $\{\beta^{k}\}$ such that $\lim\limits_{k(\in\mathcal{K})\rightarrow\infty}\beta^{k}=\beta^{\infty}$ and \\ $\lim\limits_{k(\in\mathcal{K})\rightarrow\infty}g(\beta^{k})=g(\beta^{\infty})$;
\item There exist a constant $K>0$ and $\varepsilon^{k+1}\in\partial g(\beta^{k+1})$ such that $\|\varepsilon^{k+1}\|\leq K\|\beta^{k+1}-\beta^{k}\|$, for $k$ sufficiently large.
\end{enumerate}
The first two properties are already known. To establish the third property, let $\bar{\varepsilon}^{k+1}=\nabla q_{2}(\beta^{k})-\nabla q_{2}(\beta^{k+1})-\tau_{2,k}(\beta^{k+1}-\beta^{k})-\iota_{2,k}X^{T}(X(\beta^{k+1}-\beta^{k})+r^{k+1})$. Then $\bar{\varepsilon}^{k+1}\in\partial g_{k}(\beta^{k+1})$ with $g_{k}(\beta)=h(X\beta-b+r^{k+1})+\lambda q_{1}(\beta)-q_{2}(\beta)$. Since the function $h$ is polyhedral, it is known from Theorem 23.8 and Theorem 23.9 of \cite{Rockafellar70} that $\partial g_{k}(\beta^{k+1})=X^{T}\partial h(X\beta^{k+1}-b+r^{k+1}) +\partial q(\beta^{k+1})$, where $q(\beta):=\lambda q_{1}(\beta)-q_{2}(\beta)$. Furthermore, $\partial h$ is a polyhedral multifunction and locally upper Lipschitz continuous at the point $X\beta^{k+1}-b$. Therefore, there exists $\kappa>0$ such that $\partial h(X\beta^{k+1}-b+r^{k+1})\subseteq\partial h(X\beta^{k+1}-b)+\kappa\|r^{k+1}\|\mathbb{B}_{p}$ with $\mathbb{B}_{p}$ denoting the unit ball in $\mathcal{R}^{p}$. We can find $d\bar{\varepsilon}^{k+1}\in\mathcal{R}^{p}$ with $\|d\bar{\varepsilon}^{k+1}\|\leq\kappa\|X\|\|r^{k+1}\|$ such that $\varepsilon^{k+1}:=\bar{\varepsilon}^{k+1}+d\bar{\varepsilon}^{k+1}\in\partial g(\beta^{k+1})$. Since $\nabla q_{2}$ is locally Lipschitz continuous at all $\beta\in\mathcal{C}^{\infty}$ and $\|r^{k+1}\|\leq\frac{1}{\sqrt{2}}\|X\|\|\beta^{k+1}-\beta^{k}\|$ according to \eqref{cond:algm-ppmm}, we can find some $K>0$ such that $\|\varepsilon^{k+1}\|\leq K\|\beta^{k+1}-\beta^{k}\|$. Combining properties (i)-(iii) and Proposition 4 of \cite{bolte-pauwels2016}, the desired result follows easily.

\end{proof}

\section{Numerical experiments}
\label{sec:Numerical experiments}
In this section, we implement some numerical experiments to demonstrate the efficiency of our PPMM algorithm for the tuning-free robust regression problems. The numerical experiments are implemented on a PC (Intel Core 2 Duo 2.6 GHz with 4 GB RAM) on two types of data sets. The first type of data set is generated randomly in the low sample high dimension setting which is similar to that in \cite{WangPeng}. The synthetic data is from the regression model
\begin{eqnarray*}
Y_{i}=X_{i}^{T}\dot{\beta}+\varsigma_{i},\ i=1,2,\ldots,n,
\end{eqnarray*}
where $X_{i}$ is generated from a $p$-dimensional multivariate normal distribution $N_{p}(0,\Sigma)$ and is independent of $\varsigma_{i}$. In the first six examples, we set $\dot{\beta}=(\sqrt{3},\sqrt{3},\sqrt{3},\textbf{0}_{p-3})$ with $\textbf{0}_{p-3}$ a $(p-3)$-dimensional vector of zeros. The covariance matrix $\Sigma$ is symmetric with $\Sigma_{ij}=0.5$ for $i\neq j$ and $\Sigma_{ij}=1$ for $i=j$. The noise $\varsigma_{i}$ is from the following different distributions.
\begin{enumerate}[(1)]
\item The normal distribution with mean 0 and variance 0.25;
\item The normal distribution with mean 0 and variance 1;
\item The normal distribution with mean 0 and variance 2;
\item The mixture normal distribution $0.95N(0,1)+0.05N(0,100)$;
\item The $t$ distribution with 4 degree of freedom;
\item The standard Cauchy distribution.
\end{enumerate}
The settings in the next six examples are the same as that in the first six examples except that $\dot{\beta}=(2,2,2,2,1.75,1.75,1.75,1.5,1.5,1.5,1.25,1.25,1.25,1,1,1,0.75,0.75,0.75,0.5,$  $0.5,0.5,0.25,0.25,0.25,\textbf{0}_{p-25})$. We set the random seed of the $i$-th example as $i$ ($i=1,2,\ldots,12$).

We also implement our experiments on some large-scale data sets obtained from the KEEL-dataset repository \cite{Alcala2011}.
As that in \cite{LiSunToh2018}, we apply the method in \cite{Huang2010} and use the polynomial basis functions to expand the features of those data sets with the last digit in the names of the data sets , e.g., baseball5, concrete7, dee10 and puma32h51tst3 denoting the order of the polynomial used to expand the features. Let $\hat{\beta}$ be the vector obtained by sorting $\beta$ such that $|\hat{\beta}_{1}|\geq|\hat{\beta}_{2}|\geq\cdots\geq|\hat{\beta}_{p}|$. In our numerical experiments, we define the number of nonzero elements as the minimal $k$ satisfying
\begin{eqnarray*}
\sum_{i=1}^{k}|\hat{\beta}_{i}|\geq 0.9999\|\beta\|_{1}.
\end{eqnarray*}
For simplicity, we use ``$s\ \mbox{sign}(t)|t|$'' to denote
a number of the form ``$s\times 10^{t}$'', e.g., 1.0-3 denotes $1.0\times 10^{-3}$.

\subsection{Numerical experiments for the $\ell_{1}$ regularized tuning-free robust regression problem}
In the subsection, we implement some performances of the ADMM, Gurobi and our proposed PPMM algorithm of Stage 1 for the $\ell_{1}$ regularized tuning-free robust regression problem (\ref{Wilcoxon-score-loss}). We adopt the relative KKT residual
\begin{eqnarray*}
\eta_{kkt}:=\max\left\{\frac{\|y-\mbox{Prox}_{h}(u+y)\|}{1+\|y\|},\frac{\|x-\mbox{Prox}_{\lambda\|\cdot\|_{1}}(x-A^{T}u)\|}{1+\|x\|},\frac{\|Ax-y-b\|}{1+\|y\|}\right\}
\end{eqnarray*}
to measure the accuracy of our PPMM algorithm and the ADMM. Our PPMM algorithm and the ADMM are terminated if the desired relative KKT residual $\eta_{kkt}<\mbox{Tol}:=10^{-6}$ or the number of the iterations reaches the maximum of $N_{\max}$. In our performances, we set $N_{max}$ as 200 for the PPMM algorithm and 40000 for the ADMM. As for Gurobi, we use the barrier algorithm (without presolve and crossover) and set the relevant tolerance BarConvTol$=10^{-6}$. All of these performances are stopped if the running time reaches the pre-set maximum of 4 hours. The tuning parameter $\lambda$ is obtained by simulation based on 1000 repetitions with the same parameter setting as that in \cite{WangPeng}.

\subsubsection{ADMM}
In this subsection, we introduce the implementation details of the ADMM for the problem (\ref{Wilcoxon-score-loss}). We first introduce auxiliary variables and reformulate (\ref{Wilcoxon-score-loss}) as follows
\begin{eqnarray}\label{admm-problem}
\min_{\beta,z\in\mathcal{R}^{p},y\in\mathcal{R}^{n}}\Big\{h(y)+\lambda\|z\|_{1}\ \Big|\ X\beta-y-b=0,\ \beta-z=0\Big\}.
\end{eqnarray}
Given $\sigma>0$, the augmented Lagrangian function corresponding to the problem (\ref{admm-problem}) takes the following form
\begin{eqnarray*}
\mathcal{L}_{\sigma}(\beta,y,z;u,v)&:=&h(y)+\lambda\|z\|_{1}+\frac{\sigma}{2}\|X\beta-y-b+\sigma^{-1}u\|^{2}-\frac{2}{\sigma}\|u\|^{2}\\
&&+\frac{\sigma}{2}\|\beta-z+\sigma^{-1}v\|-\frac{2}{\sigma}\|v\|^{2}.
\end{eqnarray*}
Therefore, we describe the details of the ADMM below.

\bigskip

\ni\fbox{\parbox{\textwidth}{\noindent{\bf Algorithm 4 (ADMM for the problem (\ref{Wilcoxon-score-loss})): \label{algorithm-ADMM}} Given
$\rho\in(0,\frac{1+\sqrt{5}}{2})$, $\sigma>0$, choose $(y^{0},z^{0},u^{0},v^{0})\in\mathcal{R}^{n}\times\mathcal{R}^{p}\times\mathcal{R}^{n}\times\mathcal{R}^{p}$. Set $k=0$ and iterate:
\begin{description}
\item [Step 1.] Compute
\begin{eqnarray*}
\beta^{k+1}&=&\mathop{\mbox{argmin}}_{\beta\in\mathcal{R}^{p}}\Big\{\mathcal{L}_{\sigma}(\beta,y^{k},z^{k};u^{k},v^{k})\Big\}\\
&=&\left(I_{p}+X^{T}X\right)^{-1}\left(z^{k}-\sigma^{-1}v^{k}+X^{T}(y^{k}+b-\sigma^{-1}u^{k})\right);\\
(y^{k+1},z^{k+1})&=&\mathop{\mbox{argmin}}_{y\in\mathcal{R}^{n},z\in\mathcal{R}^{p}}\Big\{\mathcal{L}_{\sigma}(\beta^{k+1},y,z;u^{k},v^{k})\Big\}\\
&=&\left(\mbox{Prox}_{\sigma^{-1}h}(X\beta^{k+1}-b+\sigma^{-1}u^{k}),\mbox{Prox}_{\sigma^{-1}\lambda\|\cdot\|_{1}}(\beta^{k+1}+\sigma^{-1}v^{k})\right).
\end{eqnarray*}
\item [Step 2.] Update the variables $u$ and $v$.
\begin{eqnarray*}
u^{k+1}=u^{k}+\rho\sigma(X\beta^{k+1}-y^{k+1}-b),\quad v^{k+1}=v^{k}+\rho\sigma(\beta^{k+1}-z^{k+1}).
\end{eqnarray*}
If the corresponding stopping criteria is satisfied, then terminate; otherwise set $k:=k+1$ and goto Step 1.
\end{description}}}

\bigskip

\subsubsection{Solving the problem (\ref{Wilcoxon-score-loss}) by linear programming}
The problem (\ref{Wilcoxon-score-loss}) is polyhedral and can be solved via linear programming. By introducing auxiliary variables, the problem (\ref{Wilcoxon-score-loss}) can be equivalently written as the following form
\begin{eqnarray*}
&&\min_{\beta,\xi,\zeta}\left\{\frac{1}{n(n-1)}\sum_{1\leq i<j\leq n}(\xi^{+}_{ij}+\xi^{-}_{ij})+\lambda\sum_{k=1}^{p}\zeta_{k}\right\}\\
&&s.t.\ \xi^{+}_{ij}-\xi^{-}_{ij}=(b_{i}-X_{i}\beta)-(b_{j}-X_{j}\beta),\ 1\leq i<j\leq n,\\
&&\quad\ \ \xi^{+}_{ij}\geq 0,\ \xi^{-}_{ij}\geq 0,\ 1\leq i<j\leq n,\\
&&\quad\ \ \zeta_{k}\geq\beta_{k},\ \zeta_{k}\geq-\beta_{k},\ k=1,2\ldots,p.
\end{eqnarray*}
The above is a linear programming problem and can be solved by some existing optimization software packages, e.g. Gurobi. In our comparison, we use the Gurobi package to solve the above problem.


\subsubsection{Numerical results for the tuning-free robust Lasso problem (\ref{Wilcoxon-score-loss})} \label{subsubsec:convex-numerical-exp}
In this subsection, we show some performances of numerical results for the tuning-free robust Lasso problem (\ref{Wilcoxon-score-loss}). In our comparison, we report the problem (pbname), the number of samples ($n$), features ($p$) and nonzero elements (nnz), lambda ($\lambda$), the relative KKT residual ($\eta_{kkt}$), the primal objective value (pobj) and the running time (time) in the format of ``hours:minutes:seconds''. In the synthetic data sets, for a solution $\widehat{\beta}$ obtained by an algorithm, we also report the $L_{1}$ estimation error $\|\widehat{\beta}-\dot{\beta}\|_{1}$ ($L_{1}$), the $L_{2}$ estimation error $\|\widehat{\beta}-\dot{\beta}\|_{2}$ ($L_{2}$), the model error $(\widehat{\beta}-\dot{\beta})^{T}\Sigma_{X}(\widehat{\beta}-\dot{\beta})$ (ME) with $\Sigma_{X}$ the covariance matrix of $X$, the number of false positive variables (FP) which is the number of the noise covariates that are selected in the model and the number of false negative variables (FN) which is the number of the active variables that are not selected in the model. Except the three columns of results including $\eta_{kkt}$, pobj and time, the other columns of results are obtained by the PPMM algorithm.

Due to the excessive memory requirement of Gurobi, we first compare the PPMM algorithm, ADMM and Gurobi on the synthetic datasets with $n=200$ and $p=800$. The results are listed in Table \ref{table-compare-random-convex-all}.

\begin{table}[!htbp]
\scriptsize
\begin{center}
\setlength{\belowcaptionskip}{10pt}
\parbox{.85\textwidth}{\caption{The performances of the PPMM, ADMM and Gurobi on the synthetic datasets ($n=200$, $p=1000$) for the tuning-free robust Lasso problem. In this table, ``$a$''=PPMM, ``$b$"=ADMM, ``$c$''=Gurobi.}\label{table-compare-random-convex-all}}
\begin{tabular}{|p{0.55cm}<{\centering}|p{0.5cm}<{\centering}|p{0.25cm}<{\centering}|p{2.6cm}<{\centering}|p{4.1cm}<{\centering}|p{1.7cm}<{\centering}|p{0.55cm}<{\centering}|p{0.55cm}<{\centering}|p{0.55cm}<{\centering}|p{0.25cm}<{\centering}|p{0.25cm}<{\centering}|}
\hline
\multirow{2}{*}{pbname}&\multirow{2}{*}{$\lambda$}&\multirow{2}{*}{nnz}&\multicolumn{1}{c|}{$\eta_{kkt}$}&\multicolumn{1}{c|}{pobj}&\multicolumn{1}{c|}{time}&\multirow{2}{*}{$L_{1}$}&\multirow{2}{*}{$L_{2}$}&\multirow{2}{*}{ME}&\multirow{2}{*}{FP}&\multirow{2}{*}{FN}\\\cline{4-6} $n$; $p$ & & & $a\ |\ b\ |\ c$ & $a\ |\ b\ |\ c$ & $a\ |\ b\ |\ c$ & & & & &\\[2pt] 
\hline\mc{1}{|c|}{exmp1} 
	 &\multirow{2}{*}{0.271}	 &\multirow{2}{*}{4}&\multirow{2}{*}{ 2.3-7 $|$  1.0-6 $|$  2.2-10} &\multirow{2}{*}{ 1.4990+0 $|$      1.4990+0  $|$      1.4990+0} &\multirow{2}{*}{02 $|$ 1:27  $|$ 1:14}  &\multirow{2}{*}{2.7-1} &\multirow{2}{*}{1.7-1} &\multirow{2}{*}{3.2-2} &\multirow{2}{*}{1} &\multirow{2}{*}{0}
\\[2pt] 
\mc{1}{|c|}{200;1000} & & & & & & & & & &\\[2pt] 
\hline\mc{1}{|c|}{exmp2} 
	 &\multirow{2}{*}{0.282}	 &\multirow{2}{*}{6}&\multirow{2}{*}{ 4.8-7 $|$  \emph{2.0-6} $|$  3.3-10} &\multirow{2}{*}{ 1.9116+0 $|$      1.9117+0  $|$      1.9116+0} &\multirow{2}{*}{02 $|$ 1:54  $|$ 45}  &\multirow{2}{*}{1.4+0} &\multirow{2}{*}{7.1-1} &\multirow{2}{*}{5.8-1} &\multirow{2}{*}{3} &\multirow{2}{*}{0}
\\[2pt] 
\mc{1}{|c|}{200;1000} & & & & & & & & & &\\[2pt] 
\hline\mc{1}{|c|}{exmp3} 
	 &\multirow{2}{*}{0.267}	 &\multirow{2}{*}{5}&\multirow{2}{*}{ 9.8-7 $|$  1.0-6 $|$  8.5-10} &\multirow{2}{*}{ 2.2355+0 $|$      2.2355+0  $|$      2.2355+0} &\multirow{2}{*}{02 $|$ 52  $|$ 53}  &\multirow{2}{*}{2.0+0} &\multirow{2}{*}{1.1+0} &\multirow{2}{*}{2.0+0} &\multirow{2}{*}{2} &\multirow{2}{*}{0}
\\[2pt] 
\mc{1}{|c|}{200;1000} & & & & & & & & & &\\[2pt] 
\hline\mc{1}{|c|}{exmp4} 
	 &\multirow{2}{*}{0.272}	 &\multirow{2}{*}{4}&\multirow{2}{*}{ 3.6-7 $|$  \emph{1.3-4} $|$  1.1-9} &\multirow{2}{*}{ 3.6425+0 $|$      3.6437+0  $|$      3.6425+0} &\multirow{2}{*}{02 $|$ 2:01  $|$ 53}  &\multirow{2}{*}{4.6+0} &\multirow{2}{*}{2.5+0} &\multirow{2}{*}{1.2+1} &\multirow{2}{*}{1} &\multirow{2}{*}{0}
\\[2pt] 
\mc{1}{|c|}{200;1000} & & & & & & & & & &\\[2pt] 
\hline\mc{1}{|c|}{exmp5} 
	 &\multirow{2}{*}{0.276}	 &\multirow{2}{*}{7}&\multirow{2}{*}{ 6.3-7 $|$  \emph{3.4-6} $|$  3.4-10} &\multirow{2}{*}{ 2.1721+0 $|$      2.1721+0  $|$      2.1721+0} &\multirow{2}{*}{02 $|$ 2:17  $|$ 1:04}  &\multirow{2}{*}{1.9+0} &\multirow{2}{*}{1.0+0} &\multirow{2}{*}{1.9+0} &\multirow{2}{*}{4} &\multirow{2}{*}{0}
\\[2pt] 
\mc{1}{|c|}{200;1000} & & & & & & & & & &\\[2pt] 
\hline\mc{1}{|c|}{exmp6} 
	 &\multirow{2}{*}{0.270}	 &\multirow{2}{*}{9}&\multirow{2}{*}{ 9.9-7 $|$  \emph{1.1-5} $|$  1.4-9} &\multirow{2}{*}{ 2.9175+0 $|$      2.9177+0  $|$      2.9175+0} &\multirow{2}{*}{02 $|$ 2:26  $|$ 1:29}  &\multirow{2}{*}{2.5+0} &\multirow{2}{*}{1.2+0} &\multirow{2}{*}{2.0+0} &\multirow{2}{*}{6} &\multirow{2}{*}{0}
\\[2pt] 
\mc{1}{|c|}{200;1000} & & & & & & & & & &\\[2pt] 
\hline\mc{1}{|c|}{exmp7} 
	 &\multirow{2}{*}{0.266}	 &\multirow{2}{*}{94}&\multirow{2}{*}{ 9.9-7 $|$  \emph{1.4-5} $|$  2.0-9} &\multirow{2}{*}{ 7.7347+0 $|$      7.7352+0  $|$      7.7347+0} &\multirow{2}{*}{06 $|$ 2:42  $|$ 1:09}  &\multirow{2}{*}{1.4+1} &\multirow{2}{*}{2.0+0} &\multirow{2}{*}{6.9-1} &\multirow{2}{*}{72} &\multirow{2}{*}{3}
\\[2pt] 
\mc{1}{|c|}{200;1000} & & & & & & & & & &\\[2pt] 
\hline\mc{1}{|c|}{exmp8} 
	 &\multirow{2}{*}{0.267}	 &\multirow{2}{*}{88}&\multirow{2}{*}{ 9.5-7 $|$  \emph{5.6-5} $|$  5.8-10} &\multirow{2}{*}{ 7.8796+0 $|$      7.8816+0  $|$      7.8796+0} &\multirow{2}{*}{05 $|$ 3:00  $|$ 53}  &\multirow{2}{*}{1.9+1} &\multirow{2}{*}{2.9+0} &\multirow{2}{*}{2.3+0} &\multirow{2}{*}{68} &\multirow{2}{*}{5}
\\[2pt] 
\mc{1}{|c|}{200;1000} & & & & & & & & & &\\[2pt] 
\hline\mc{1}{|c|}{exmp9} 
	 &\multirow{2}{*}{0.268}	 &\multirow{2}{*}{79}&\multirow{2}{*}{ 4.5-7 $|$  \emph{1.5-5} $|$  3.4-9} &\multirow{2}{*}{ 8.2800+0 $|$      8.2804+0  $|$      8.2800+0} &\multirow{2}{*}{05 $|$ 3:04  $|$ 1:01}  &\multirow{2}{*}{2.4+1} &\multirow{2}{*}{3.5+0} &\multirow{2}{*}{4.4+0} &\multirow{2}{*}{59} &\multirow{2}{*}{5}
\\[2pt] 
\mc{1}{|c|}{200;1000} & & & & & & & & & &\\[2pt] 
\hline\mc{1}{|c|}{exmp10} 
	 &\multirow{2}{*}{0.269}	 &\multirow{2}{*}{59}&\multirow{2}{*}{ 5.3-7 $|$  \emph{1.1-5} $|$  9.8-10} &\multirow{2}{*}{ 9.2758+0 $|$      9.2760+0  $|$      9.2758+0} &\multirow{2}{*}{04 $|$ 3:07  $|$ 50}  &\multirow{2}{*}{3.9+1} &\multirow{2}{*}{6.3+0} &\multirow{2}{*}{2.4+1} &\multirow{2}{*}{49} &\multirow{2}{*}{15}
\\[2pt] 
\mc{1}{|c|}{200;1000} & & & & & & & & & &\\[2pt] 
\hline\mc{1}{|c|}{exmp11} 
	 &\multirow{2}{*}{0.272}	 &\multirow{2}{*}{67}&\multirow{2}{*}{ 4.0-7 $|$  \emph{2.0-5} $|$  3.4-9} &\multirow{2}{*}{ 8.5739+0 $|$      8.5747+0  $|$      8.5739+0} &\multirow{2}{*}{04 $|$ 2:56  $|$ 1:01}  &\multirow{2}{*}{1.9+1} &\multirow{2}{*}{3.0+0} &\multirow{2}{*}{3.5+0} &\multirow{2}{*}{45} &\multirow{2}{*}{3}
\\[2pt] 
\mc{1}{|c|}{200;1000} & & & & & & & & & &\\[2pt] 
\hline\mc{1}{|c|}{exmp12} 
	 &\multirow{2}{*}{0.259}	 &\multirow{2}{*}{58}&\multirow{2}{*}{ 7.8-7 $|$  \emph{4.7-5} $|$  5.2-10} &\multirow{2}{*}{ 9.0375+0 $|$      9.0383+0  $|$      9.0374+0} &\multirow{2}{*}{04 $|$ 2:52  $|$ 52}  &\multirow{2}{*}{4.0+1} &\multirow{2}{*}{6.2+0} &\multirow{2}{*}{1.7+1} &\multirow{2}{*}{50} &\multirow{2}{*}{17}
\\[2pt] 
\mc{1}{|c|}{200;1000} & & & & & & & & & &\\[2pt] 
\hline

\end{tabular}
\end{center}
\end{table}

\begin{table}[!htbp]
\scriptsize
\begin{center}
\setlength{\belowcaptionskip}{10pt}
\parbox{.85\textwidth}{\caption{The performances of the PPMM and ADMM on the synthetic datasets ($n=2000$, $p=8000$) for the tuning-free robust Lasso problem. In this table, ``$a$''=PPMM, ``$b$"=ADMM, ``$c$''=Gurobi.}\label{table-compare-random-convex}}
\begin{tabular}{|p{0.55cm}<{\centering}|p{0.5cm}<{\centering}|p{0.25cm}<{\centering}|p{2.0cm}<{\centering}|p{3.1cm}<{\centering}|p{2.0cm}<{\centering}|p{0.55cm}<{\centering}|p{0.55cm}<{\centering}|p{0.55cm}<{\centering}|p{0.3cm}<{\centering}|p{0.3cm}<{\centering}|}
\hline
\multirow{2}{*}{pbname}&\multirow{2}{*}{$\lambda$}&\multirow{2}{*}{nnz}&\multicolumn{1}{c|}{$\eta_{kkt}$}&\multicolumn{1}{c|}{pobj}&\multicolumn{1}{c|}{time}&\multirow{2}{*}{$L_{1}$}&\multirow{2}{*}{$L_{2}$}&\multirow{2}{*}{ME}&\multirow{2}{*}{FP}&\multirow{2}{*}{FN}\\\cline{4-6} $n$; $p$ & & & $a\ |\ b\ |\ c$ & $a\ |\ b\ |\ c$ & $a\ |\ b\ |\ c$ & & & & &\\[2pt] 
\hline\mc{1}{|c|}{exmp1} 
	 &\multirow{2}{*}{0.095}	 &\multirow{2}{*}{7}&\multirow{2}{*}{ 2.7-7 $|$  \emph{4.0-4} $|$ --} &\multirow{2}{*}{ 6.3448-1 $|$      6.4076-1 $|$ --} &\multirow{2}{*}{53 $|$ 1:05:46 $|$ --}  &\multirow{2}{*}{1.0-1} &\multirow{2}{*}{5.3-2} &\multirow{2}{*}{3.4-3} &\multirow{2}{*}{4} &\multirow{2}{*}{0}
\\[2pt] 
\mc{1}{|c|}{2000;8000} & & & & & & & & & &\\[2pt] 
\hline\mc{1}{|c|}{exmp2} 
	 &\multirow{2}{*}{0.094}	 &\multirow{2}{*}{14}&\multirow{2}{*}{ 3.2-7 $|$  \emph{5.7-4} $|$ --} &\multirow{2}{*}{ 1.0409+0 $|$      1.0464+0 $|$ --} &\multirow{2}{*}{50 $|$ 1:06:17 $|$ --}  &\multirow{2}{*}{5.0-1} &\multirow{2}{*}{2.3-1} &\multirow{2}{*}{6.4-2} &\multirow{2}{*}{11} &\multirow{2}{*}{0}
\\[2pt] 
\mc{1}{|c|}{2000;8000} & & & & & & & & & &\\[2pt] 
\hline\mc{1}{|c|}{exmp3} 
	 &\multirow{2}{*}{0.094}	 &\multirow{2}{*}{9}&\multirow{2}{*}{ 4.7-7 $|$  \emph{1.6-4} $|$ --} &\multirow{2}{*}{ 1.5989+0 $|$      1.6002+0 $|$ --} &\multirow{2}{*}{46 $|$ 1:06:36 $|$ --}  &\multirow{2}{*}{7.2-1} &\multirow{2}{*}{3.5-1} &\multirow{2}{*}{1.5-1} &\multirow{2}{*}{6} &\multirow{2}{*}{0}
\\[2pt] 
\mc{1}{|c|}{2000;8000} & & & & & & & & & &\\[2pt] 
\hline\mc{1}{|c|}{exmp4} 
	 &\multirow{2}{*}{0.095}	 &\multirow{2}{*}{14}&\multirow{2}{*}{ 8.6-7 $|$  \emph{1.4-3} $|$ --} &\multirow{2}{*}{ 3.3337+0 $|$      3.3628+0 $|$ --} &\multirow{2}{*}{41 $|$ 1:02:23 $|$ --}  &\multirow{2}{*}{2.1+0} &\multirow{2}{*}{9.9-1} &\multirow{2}{*}{1.4+0} &\multirow{2}{*}{11} &\multirow{2}{*}{0}
\\[2pt] 
\mc{1}{|c|}{2000;8000} & & & & & & & & & &\\[2pt] 
\hline\mc{1}{|c|}{exmp5} 
	 &\multirow{2}{*}{0.096}	 &\multirow{2}{*}{11}&\multirow{2}{*}{ 8.8-7 $|$  \emph{3.4-4} $|$ --} &\multirow{2}{*}{ 1.4533+0 $|$      1.4563+0 $|$ --} &\multirow{2}{*}{44 $|$ 1:02:22 $|$ --}  &\multirow{2}{*}{6.2-1} &\multirow{2}{*}{3.1-1} &\multirow{2}{*}{1.3-1} &\multirow{2}{*}{8} &\multirow{2}{*}{0}
\\[2pt] 
\mc{1}{|c|}{2000;8000} & & & & & & & & & &\\[2pt] 
\hline\mc{1}{|c|}{exmp6} 
	 &\multirow{2}{*}{0.093}	 &\multirow{2}{*}{13}&\multirow{2}{*}{ 9.3-7 $|$  \emph{8.9-5} $|$ --} &\multirow{2}{*}{ 7.6828+0 $|$      7.6840+0 $|$ --} &\multirow{2}{*}{44 $|$ 1:02:23 $|$ --}  &\multirow{2}{*}{8.3-1} &\multirow{2}{*}{4.2-1} &\multirow{2}{*}{2.2-1} &\multirow{2}{*}{10} &\multirow{2}{*}{0}
\\[2pt] 
\mc{1}{|c|}{2000;8000} & & & & & & & & & &\\[2pt] 
\hline\mc{1}{|c|}{exmp7} 
	 &\multirow{2}{*}{0.095}	 &\multirow{2}{*}{144}&\multirow{2}{*}{ 1.0-6 $|$  \emph{2.0-5} $|$ --} &\multirow{2}{*}{ 2.8776+0 $|$      2.8785+0 $|$ --} &\multirow{2}{*}{4:59 $|$ 1:02:20 $|$ --}  &\multirow{2}{*}{9.2-1} &\multirow{2}{*}{1.2-1} &\multirow{2}{*}{7.9-3} &\multirow{2}{*}{119} &\multirow{2}{*}{0}
\\[2pt] 
\mc{1}{|c|}{2000;8000} & & & & & & & & & &\\[2pt] 
\hline\mc{1}{|c|}{exmp8} 
	 &\multirow{2}{*}{0.093}	 &\multirow{2}{*}{147}&\multirow{2}{*}{ 3.8-7 $|$  \emph{8.4-5} $|$ --} &\multirow{2}{*}{ 3.2321+0 $|$      3.2361+0 $|$ --} &\multirow{2}{*}{3:28 $|$ 1:02:30 $|$ --}  &\multirow{2}{*}{3.4+0} &\multirow{2}{*}{4.6-1} &\multirow{2}{*}{1.4-1} &\multirow{2}{*}{122} &\multirow{2}{*}{0}
\\[2pt] 
\mc{1}{|c|}{2000;8000} & & & & & & & & & &\\[2pt] 
\hline\mc{1}{|c|}{exmp9} 
	 &\multirow{2}{*}{0.095}	 &\multirow{2}{*}{168}&\multirow{2}{*}{ 4.0-7 $|$  \emph{6.4-5} $|$ --} &\multirow{2}{*}{ 3.8115+0 $|$      3.8143+0 $|$ --} &\multirow{2}{*}{3:30 $|$ 1:02:20 $|$ --}  &\multirow{2}{*}{7.0+0} &\multirow{2}{*}{8.9-1} &\multirow{2}{*}{5.2-1} &\multirow{2}{*}{143} &\multirow{2}{*}{0}
\\[2pt] 
\mc{1}{|c|}{2000;8000} & & & & & & & & & &\\[2pt] 
\hline\mc{1}{|c|}{exmp10} 
	 &\multirow{2}{*}{0.096}	 &\multirow{2}{*}{157}&\multirow{2}{*}{ 7.3-7 $|$  \emph{1.9-4} $|$ --} &\multirow{2}{*}{ 5.5081+0 $|$      5.5155+0 $|$ --} &\multirow{2}{*}{2:49 $|$ 1:02:14 $|$ --}  &\multirow{2}{*}{1.8+1} &\multirow{2}{*}{2.3+0} &\multirow{2}{*}{3.7+0} &\multirow{2}{*}{135} &\multirow{2}{*}{3}
\\[2pt] 
\mc{1}{|c|}{2000;8000} & & & & & & & & & &\\[2pt] 
\hline\mc{1}{|c|}{exmp11} 
	 &\multirow{2}{*}{0.094}	 &\multirow{2}{*}{136}&\multirow{2}{*}{ 2.5-7 $|$  \emph{5.7-5} $|$ --} &\multirow{2}{*}{ 3.6788+0 $|$      3.6815+0 $|$ --} &\multirow{2}{*}{2:58 $|$ 1:02:02 $|$ --}  &\multirow{2}{*}{5.3+0} &\multirow{2}{*}{7.1-1} &\multirow{2}{*}{3.5-1} &\multirow{2}{*}{111} &\multirow{2}{*}{0}
\\[2pt] 
\mc{1}{|c|}{2000;8000} & & & & & & & & & &\\[2pt] 
\hline\mc{1}{|c|}{exmp12} 
	 &\multirow{2}{*}{0.093}	 &\multirow{2}{*}{163}&\multirow{2}{*}{ 9.4-7 $|$  \emph{2.5-4} $|$ --} &\multirow{2}{*}{ 1.4266+1 $|$      1.4276+1 $|$ --} &\multirow{2}{*}{2:50 $|$ 1:02:15 $|$ --}  &\multirow{2}{*}{8.2+0} &\multirow{2}{*}{1.1+0} &\multirow{2}{*}{7.6-1} &\multirow{2}{*}{138} &\multirow{2}{*}{0}
\\[2pt] 
\mc{1}{|c|}{2000;8000} & & & & & & & & & &\\[2pt] 
\hline

\end{tabular}
\end{center}
\end{table}

\begin{table}[!htbp]
\scriptsize
\begin{center}
\setlength{\belowcaptionskip}{10pt}
\parbox{0.75\textwidth}{\caption{The performances of the PPMM and ADMM on the KEEL datasets for the tuning-free robust Lasso problem. In the table, ``$a$''=PPMM, ``$b$"=ADMM, ``$c$''=Gurobi.}\label{table-compare-KEEL-convex}}
\begin{tabular}{|p{2.5cm}<{\centering}|p{0.6cm}<{\centering}|p{0.4cm}<{\centering}|p{2.4cm}<{\centering}|p{3.2cm}<{\centering}|p{2.7cm}<{\centering}|}
\hline
\multirow{2}{*}{pbname}&\multirow{2}{*}{$\lambda$}&\multirow{2}{*}{nnz}&\multicolumn{1}{c|}{$\eta_{kkt}$}&\multicolumn{1}{c|}{pobj}&\multicolumn{1}{c|}{time}\\\cline{4-6} $n$; $p$ & & & $a\ |\ b\ |\ c$ & $a\ |\ b\ |\ c$ & $a\ |\ b\ |\ c$\\[2pt] 
\hline\mc{1}{|c|}{baseball5} 
	 &\multirow{2}{*}{0.141}	 &\multirow{2}{*}{38}&\multirow{2}{*}{ 1.0-6 $|$  \emph{3.8-5} $|$ --} &\multirow{2}{*}{ 5.8139+2 $|$      5.6807+2 $|$ --} &\multirow{2}{*}{16 $|$ 36:25 $|$ --}
\\[2pt] 
\mc{1}{|c|}{337;20349} & & & & & \\[2pt] 
\hline\mc{1}{|c|}{compactiv51tra4} 
	 &\multirow{2}{*}{0.014}	 &\multirow{2}{*}{36}&\multirow{2}{*}{ 8.6-7 $|$  \emph{2.1-2} $|$ --} &\multirow{2}{*}{ 3.3341+0 $|$      5.4468+0 $|$ --} &\multirow{2}{*}{27:58 $|$ 4:01:28 $|$ --}
\\[2pt] 
\mc{1}{|c|}{6553;12650} & & & & & \\[2pt] 
\hline\mc{1}{|c|}{concrete7} 
	 &\multirow{2}{*}{0.046}	 &\multirow{2}{*}{16}&\multirow{2}{*}{ 9.6-7 $|$  1.0-6 $|$ --} &\multirow{2}{*}{ 7.1762+0 $|$      7.1719+0 $|$ --} &\multirow{2}{*}{48 $|$ 7:35 $|$ --}
\\[2pt] 
\mc{1}{|c|}{1030;6435} & & & & & \\[2pt] 
\hline\mc{1}{|c|}{dee10} 
	 &\multirow{2}{*}{0.060}	 &\multirow{2}{*}{5}&\multirow{2}{*}{ 9.2-7 $|$  9.7-7 $|$ --} &\multirow{2}{*}{ 3.4835-1 $|$      3.4836-1 $|$ --} &\multirow{2}{*}{11 $|$ 1:28 $|$ --}
\\[2pt] 
\mc{1}{|c|}{365;8008} & & & & & \\[2pt] 
\hline\mc{1}{|c|}{friedman10} 
	 &\multirow{2}{*}{0.039}	 &\multirow{2}{*}{7}&\multirow{2}{*}{ 5.8-7 $|$  9.7-7 $|$ --} &\multirow{2}{*}{ 1.8615+0 $|$      1.8615+0 $|$ --} &\multirow{2}{*}{16 $|$ 1:40 $|$ --}
\\[2pt] 
\mc{1}{|c|}{1200;3003} & & & & & \\[2pt] 
\hline\mc{1}{|c|}{mortgage5} 
	 &\multirow{2}{*}{0.040}	 &\multirow{2}{*}{9}&\multirow{2}{*}{ 4.6-7 $|$  1.0-6 $|$ --} &\multirow{2}{*}{ 3.8343-1 $|$      3.8345-1 $|$ --} &\multirow{2}{*}{2:34 $|$ 38:03 $|$ --}
\\[2pt] 
\mc{1}{|c|}{1049;15504} & & & & & \\[2pt] 
\hline\mc{1}{|c|}{puma32h51tst3} 
	 &\multirow{2}{*}{0.046}	 &\multirow{2}{*}{3}&\multirow{2}{*}{ 7.0-7 $|$  \emph{6.9-5} $|$ --} &\multirow{2}{*}{ 1.3671-2 $|$      1.3792-2 $|$ --} &\multirow{2}{*}{18 $|$ 58:22 $|$ --}
\\[2pt] 
\mc{1}{|c|}{1639;6545} & & & & & \\[2pt] 
\hline\mc{1}{|c|}{wizmir7} 
	 &\multirow{2}{*}{0.027}	 &\multirow{2}{*}{8}&\multirow{2}{*}{ 8.7-7 $|$  1.0-6 $|$ --} &\multirow{2}{*}{ 1.6023+0 $|$      1.5853+0 $|$ --} &\multirow{2}{*}{42 $|$ 16:03 $|$ --}
\\[2pt] 
\mc{1}{|c|}{1461;11440} & & & & & \\[2pt] 
\hline

\end{tabular}
\end{center}
\end{table}

In the subsequent large-scale experiments,
we list the performances of the comparisons in Tables \ref{table-compare-random-convex} and \ref{table-compare-KEEL-convex}. The symbol ``--" in Tables \ref{table-compare-random-convex} and \ref{table-compare-KEEL-convex} means that Gurobi fails to solve the problem due to excessive memory requirement. In the implementation of the ADMM, we apply the Sherman-Morrison-Woodbury formula \cite{Golub1996} if it is necessary, depending on the size of $n$ and $p$. We solve the linear system either by the Cholesky factorization or by an iterative solver such as the preconditioned conjugate gradient (PCG) method.  From the comparisons, we can see that the PPMM algorithm can obtain all the solutions with a desired accuracy efficiently. When the sample size $n$ and the number of features $p$ are small, Gurobi can also solve the problem with high accuracy but with much more time. For high dimensional problems, Gurobi fails with too much memory consuming.

\subsection{Numerical experiments for the nonconvex tuning-free robust regression problem}
In this section, we implement some experiments for the tuning-free robust regression problems with nonconvex regularizers. As we can see from Section \ref{subsubsec:convex-numerical-exp} that Gurobi fails to solve high dimensional problems and it is time consuming. Although we can reformulate each convex subproblem as a linear programming problem and solve it by Gurobi, we only compare the the performances of the ADMM and our PPMM algorithm for the nonconvex tuning-free robust regression problem (\ref{tuning-free robust-problem}). We adopt the relative KKT residual
\begin{eqnarray*}
\hat{\eta}_{kkt}:=\max\left\{\frac{\|y-\mbox{Prox}_{h}(u+y)\|}{1+\|y\|},\frac{\|x-\mbox{Prox}_{\lambda q_{1}-q_{2}}(x-A^{T}u)\|}{1+\|x\|},\frac{\|Ax-y-b\|}{1+\|y\|}\right\}
\end{eqnarray*}
to measure the accuracy of the performance.

We first implement Stage 1 of our PPMM algorithm with $\eta_{kkt}<10^{-4}$ to generate an initial point for Stage 2. All of the performances are stopped if the relative KKT residual $\hat{\eta}_{kkt}<10^{-6}$.  In addition, all of the algorithms are also stopped when they reach the pre-set maximum number of iterations (200 for Stage 2 of the
PPMM algorithm and 40000 for the ADMM) or the pre-set maximum running time of 4 hours. For each
synthetic data set, we generate 2000 observations for the training data set and 400 observations for the validation data set. By applying the PPMM algorithm, we fit the model on the training data set and use the validation data set to select the tuning parameter $\lambda$. For each KEEL data set, we use the PPMM algorithm and adopt a tenfold
cross validation to find the parameter $\lambda$.

\subsubsection{The ADMM for the nonconvex tuning-free robust regression problem}
In this subsection, we describe the ADMM for solving the nonconvex tuning-free robust regression problem. By introducing auxiliary variables, we can reformulate the problem (\ref{tuning-free robust-problem}) as below
\begin{eqnarray}\label{tuning-free robust-problem-reformulate}
\min_{\beta,z\in\mathcal{R}^{p},y\in\mathcal{R}^{n}}\Big\{h(y)+\lambda q_{1}(z)-q_{2}(z)\ \Big|\ X\beta-y-b=0,\ \beta-z=0\Big\}.
\end{eqnarray}
Given $\sigma>0$ and $(u,v)\in\mathcal{R}^{n}\times\mathcal{R}^{p}$, the augmented Lagrangian function related to (\ref{tuning-free robust-problem-reformulate}) is defined by
\begin{eqnarray*}
\widehat{\mathcal{L}}_{\sigma}(\beta,y,z;u,v)&:=&h(y)+\lambda q_{1}(z)-q_{2}(z)+\frac{\sigma}{2}\|X\beta-y-b+\sigma^{-1}u\|^{2}-\frac{1}{2\sigma}\|u\|^{2}\\
&&+\frac{\sigma}{2}\|\beta-z+\sigma^{-1}v\|^{2}-\frac{1}{2\sigma}\|v\|^{2}.
\end{eqnarray*}

We list the performance of the ADMM which is not guaranteed to converge due to the nonconvexity in the following form.

\bigskip

\ni\fbox{\parbox{\textwidth}{\noindent{\bf Algorithm 5 (ADMM for the problem (\ref{tuning-free robust-problem})): \label{algorithm-ADMM-nonconvex}} Given $\sigma>0$, choose $(y^{0},z^{0},u^{0},v^{0})\in\mathcal{R}^{n}\times\mathcal{R}^{p}\times\mathcal{R}^{n}\times\mathcal{R}^{p}$. Set $k=0$ and iterate:
\begin{description}
\item [Step 1.] Compute
\begin{eqnarray*}
\beta^{k+1}&=&\mathop{\mbox{argmin}}_{\beta\in\mathcal{R}^{p}}\Big\{\widehat{\mathcal{L}}_{\sigma}(\beta,y^{k},z^{k};u^{k},v^{k})\Big\}\\
&=&\left(I_{p}+X^{T}X\right)^{-1}\left(z^{k}-\sigma^{-1}v^{k}+X^{T}(y^{k}+b-\sigma^{-1}u^{k})\right);\\
(y^{k+1},z^{k+1})&=&\mathop{\mbox{argmin}}_{y\in\mathcal{R}^{n},z\in\mathcal{R}^{p}}\Big\{\widehat{\mathcal{L}}_{\sigma}(\beta^{k+1},y,z;u^{k},v^{k})\Big\}\\
&=&\left(\mbox{Prox}_{\sigma^{-1}h}(X\beta^{k+1}-b+\sigma^{-1}u^{k}),\mbox{Prox}_{\sigma^{-1}(\lambda q_{1}-q_{2})}(\beta^{k+1}+\sigma^{-1}v^{k})\right).
\end{eqnarray*}
\item [Step 2.] Update the variables $u$ and $v$.
\begin{eqnarray*}
u^{k+1}=u^{k}+\sigma(X\beta^{k+1}-y^{k+1}-b),\quad v^{k+1}=v^{k}+\sigma(\beta^{k+1}-z^{k+1}).
\end{eqnarray*}
If the corresponding stopping criterion is satisfied, then terminate; otherwise set $k:=k+1$ and return to Step 1.
\end{description}}}

\bigskip

\subsubsection{Numerical results}
We first implement the numerical results for the tuning-free robust regression problem with the SCAD regularization. The comparison results are listed in Tables \ref{table-compare-random-SCAD} and \ref{table-compare-SCAD-KEEL}. Comparing with the $\ell_{1}$ regularizer, the SCAD regularizer achieves a sparse estimation with a fewer number of nonzeros for the synthetic type of data set. The corresponding solution is much closer to $\dot{\beta}$ with less values of $L_{1}$, $L_{2}$, ME and FP. For the first 8 synthetic data sets, the SCAD regularizer can recover the true nonzero positions. The ADMM can solve some problems with the desired accuracy, but it takes much more time than the PPMM algorithm. From the comparison, we can see that the PPMM algorithm can solve all the listed problems efficiently.

\begin{table}[!htbp]
\begin{center}
\scriptsize
\setlength{\belowcaptionskip}{10pt}
\parbox{.7\textwidth}{\caption{The performances of the PPMM and ADMM on the synthetic datasets for the SCAD regularization. In this table,``$a$"=PPMM, ``$b$''=ADMM.}\label{table-compare-random-SCAD}}
\begin{tabular}{|p{0.8cm}<{\centering}|p{0.5cm}<{\centering}|p{0.3cm}<{\centering}|p{1.6cm}<{\centering}|p{2.6cm}<{\centering}|p{1.4cm}<{\centering}|p{0.55cm}<{\centering}|p{0.55cm}<{\centering}|p{0.55cm}<{\centering}|p{0.3cm}<{\centering}|p{0.3cm}<{\centering}|}
\hline
\multirow{2}{*}{pbname}&\multirow{2}{*}{$\lambda$}&\multirow{2}{*}{nnz}&\multicolumn{1}{c|}{$\eta_{kkt}$}&\multicolumn{1}{c|}{pobj}&\multicolumn{1}{c|}{time}&\multirow{2}{*}{$L_{1}$}&\multirow{2}{*}{$L_{2}$}&\multirow{2}{*}{ME}&\multirow{2}{*}{FP}&\multirow{2}{*}{FN}\\\cline{4-6} $n$; $p$ & & & $a\ |\ b $ & $a\ |\ b $ & $a\ |\ b $ & & & & &\\[2pt] 
\hline\mc{1}{|c|}{exmp1} 
	 &\multirow{2}{*}{0.329}	 &\multirow{2}{*}{3}&\multirow{2}{*}{ 9.0-8 $|$  1.0-6} &\multirow{2}{*}{ 9.0670-1 $|$ 9.0679-1} &\multirow{2}{*}{29 $|$ 34:42} & \multirow{2}{*}{3.4-2} & \multirow{2}{*}{2.0-2} & \multirow{2}{*}{2.1-4} & \multirow{2}{*}{0} &\multirow{2}{*}{0}
\\[2pt] 
\mc{1}{|c|}{2000;8000} & & & & & & & & & &\\[2pt] 
\hline\mc{1}{|c|}{exmp2} 
	 &\multirow{2}{*}{0.302}	 &\multirow{2}{*}{3}&\multirow{2}{*}{ 3.4-7 $|$  \emph{1.8-6}} &\multirow{2}{*}{ 1.2112+0 $|$ 1.2113+0} &\multirow{2}{*}{26 $|$ 1:01:01} & \multirow{2}{*}{5.4-2} & \multirow{2}{*}{3.4-2} & \multirow{2}{*}{2.1-3} & \multirow{2}{*}{0} &\multirow{2}{*}{0}
\\[2pt] 
\mc{1}{|c|}{2000;8000} & & & & & & & & & &\\[2pt] 
\hline\mc{1}{|c|}{exmp3} 
	 &\multirow{2}{*}{0.149}	 &\multirow{2}{*}{3}&\multirow{2}{*}{ 7.5-7 $|$  \emph{4.8-6}} &\multirow{2}{*}{ 1.2830+0 $|$ 1.2832+0} &\multirow{2}{*}{24 $|$ 1:02:17} & \multirow{2}{*}{2.2-1} & \multirow{2}{*}{1.3-1} & \multirow{2}{*}{1.0-2} & \multirow{2}{*}{0} &\multirow{2}{*}{0}
\\[2pt] 
\mc{1}{|c|}{2000;8000} & & & & & & & & & &\\[2pt] 
\hline\mc{1}{|c|}{exmp4} 
	 &\multirow{2}{*}{0.131}	 &\multirow{2}{*}{3}&\multirow{2}{*}{ 5.6-7 $|$  \emph{5.3-5}} &\multirow{2}{*}{ 3.0223+0 $|$ 3.0241+0} &\multirow{2}{*}{17 $|$ 1:04:18} & \multirow{2}{*}{1.5-1} & \multirow{2}{*}{1.1-1} & \multirow{2}{*}{8.7-3} & \multirow{2}{*}{0} &\multirow{2}{*}{0}
\\[2pt] 
\mc{1}{|c|}{2000;8000} & & & & & & & & & &\\[2pt] 
\hline\mc{1}{|c|}{exmp5} 
	 &\multirow{2}{*}{0.036}	 &\multirow{2}{*}{3}&\multirow{2}{*}{ 6.4-7 $|$  \emph{4.9-1}} &\multirow{2}{*}{ 9.8737-1 $|$ 6.1499+0} &\multirow{2}{*}{30 $|$ 1:04:38} & \multirow{2}{*}{8.9-2} & \multirow{2}{*}{5.3-2} & \multirow{2}{*}{1.6-3} & \multirow{2}{*}{0} &\multirow{2}{*}{0}
\\[2pt] 
\mc{1}{|c|}{2000;8000} & & & & & & & & & &\\[2pt] 
\hline\mc{1}{|c|}{exmp6} 
	 &\multirow{2}{*}{0.171}	 &\multirow{2}{*}{3}&\multirow{2}{*}{ 5.6-7 $|$  \emph{6.8-1}} &\multirow{2}{*}{ 7.4302+0 $|$ 7.3998+1} &\multirow{2}{*}{23 $|$ 1:02:58} & \multirow{2}{*}{1.6-1} & \multirow{2}{*}{1.2-1} & \multirow{2}{*}{8.6-3} & \multirow{2}{*}{0} &\multirow{2}{*}{0}
\\[2pt] 
\mc{1}{|c|}{2000;8000} & & & & & & & & & &\\[2pt] 
\hline\mc{1}{|c|}{exmp7} 
	 &\multirow{2}{*}{0.053}	 &\multirow{2}{*}{25}&\multirow{2}{*}{ 3.5-8 $|$  \emph{1.7-2}} &\multirow{2}{*}{ 3.0300-1 $|$ 7.6740-1} &\multirow{2}{*}{58 $|$ 1:04:54} & \multirow{2}{*}{1.5-1} & \multirow{2}{*}{3.6-2} & \multirow{2}{*}{6.3-4} & \multirow{2}{*}{0} &\multirow{2}{*}{0}
\\[2pt] 
\mc{1}{|c|}{2000;8000} & & & & & & & & & &\\[2pt] 
\hline\mc{1}{|c|}{exmp8} 
	 &\multirow{2}{*}{0.040}	 &\multirow{2}{*}{25}&\multirow{2}{*}{ 1.4-7 $|$  \emph{1.5-1}} &\multirow{2}{*}{ 6.5610-1 $|$ 2.7542+0} &\multirow{2}{*}{44 $|$ 1:04:09} & \multirow{2}{*}{6.4-1} & \multirow{2}{*}{1.6-1} & \multirow{2}{*}{1.3-2} & \multirow{2}{*}{0} &\multirow{2}{*}{0}
\\[2pt] 
\mc{1}{|c|}{2000;8000} & & & & & & & & & &\\[2pt] 
\hline\mc{1}{|c|}{exmp9} 
	 &\multirow{2}{*}{0.032}	 &\multirow{2}{*}{27}&\multirow{2}{*}{ 7.4-7 $|$  \emph{3.8-1}} &\multirow{2}{*}{ 1.1848+0 $|$ 5.0077+0} &\multirow{2}{*}{38 $|$ 1:02:15} & \multirow{2}{*}{1.9+0} & \multirow{2}{*}{4.8-1} & \multirow{2}{*}{1.1-1} & \multirow{2}{*}{3} &\multirow{2}{*}{1}
\\[2pt] 
\mc{1}{|c|}{2000;8000} & & & & & & & & & &\\[2pt] 
\hline\mc{1}{|c|}{exmp10} 
	 &\multirow{2}{*}{0.032}	 &\multirow{2}{*}{24}&\multirow{2}{*}{ 9.2-7 $|$  \emph{5.5-1}} &\multirow{2}{*}{ 2.9424+0 $|$ 5.6558+0} &\multirow{2}{*}{25 $|$ 1:03:51} & \multirow{2}{*}{7.1+0} & \multirow{2}{*}{1.7+0} & \multirow{2}{*}{1.4+0} & \multirow{2}{*}{4} &\multirow{2}{*}{5}
\\[2pt] 
\mc{1}{|c|}{2000;8000} & & & & & & & & & &\\[2pt] 
\hline\mc{1}{|c|}{exmp11} 
	 &\multirow{2}{*}{0.040}	 &\multirow{2}{*}{25}&\multirow{2}{*}{ 9.0-7 $|$  \emph{2.2-1}} &\multirow{2}{*}{ 1.0917+0 $|$ 4.2131+0} &\multirow{2}{*}{33 $|$ 1:03:51} & \multirow{2}{*}{1.6+0} & \multirow{2}{*}{4.4-1} & \multirow{2}{*}{9.3-2} & \multirow{2}{*}{1} &\multirow{2}{*}{1}
\\[2pt] 
\mc{1}{|c|}{2000;8000} & & & & & & & & & &\\[2pt] 
\hline\mc{1}{|c|}{exmp12} 
	 &\multirow{2}{*}{0.081}	 &\multirow{2}{*}{20}&\multirow{2}{*}{ 2.2-7 $|$  \emph{6.9-1}} &\multirow{2}{*}{ 1.1994+1 $|$ 1.1195+2} &\multirow{2}{*}{31 $|$ 1:08:41} & \multirow{2}{*}{4.4+0} & \multirow{2}{*}{1.2+0} & \multirow{2}{*}{6.5-1} & \multirow{2}{*}{0} &\multirow{2}{*}{5}
\\[2pt] 
\mc{1}{|c|}{2000;8000} & & & & & & & & & &\\[2pt] 
\hline

\end{tabular}
\end{center}
\end{table}

\newpage
\begin{table}[!htbp]
\begin{center}
\scriptsize
\setlength{\belowcaptionskip}{10pt}
\parbox{.7\textwidth}{\caption{The performances of the PPMM and ADMM on the KEEL datasets for the SCAD regularization. In this table,``$a$"=PPMM, ``$b$''=ADMM.}\label{table-compare-SCAD-KEEL}}
\begin{tabular}{|c|c|c|c|c|c|c|}
\hline
\multirow{2}{*}{probname}&\multirow{2}{*}{$\lambda$}&\multirow{2}{*}{nnz}&\multicolumn{1}{c|}{$\eta_{kkt}$}&\multicolumn{1}{c|}{pobj}&\multicolumn{1}{c|}{time}\\\cline{4-6} $n$; $p$ & & & $a\ |\ b$ & $a\ |\ b$ & $a\ |\ b$\\[2pt] 
\hline\mc{1}{|c|}{baseball5} 
	 &\multirow{2}{*}{0.081}	 &\multirow{2}{*}{30}&\multirow{2}{*}{ 9.7-7 $|$  \emph{7.1-1}} &\multirow{2}{*}{ 3.6945+2 $|$      3.1989+2} &\multirow{2}{*}{17 $|$ 35:53}
\\[2pt] 
\mc{1}{|c|}{337;20349} & & & & & \\[2pt] 
\hline\mc{1}{|c|}{compactiv51tra4} 
	 &\multirow{2}{*}{0.057}	 &\multirow{2}{*}{3}&\multirow{2}{*}{ 5.6-7 $|$  \emph{7.3-1}} &\multirow{2}{*}{ 4.6881+0 $|$      3.4133+1} &\multirow{2}{*}{3:38 $|$ 4:01:27}
\\[2pt] 
\mc{1}{|c|}{6553;12650} & & & & & \\[2pt] 
\hline\mc{1}{|c|}{concrete7} 
	 &\multirow{2}{*}{0.040}	 &\multirow{2}{*}{9}&\multirow{2}{*}{ 9.2-7 $|$  \emph{6.7-1}} &\multirow{2}{*}{ 4.1622+0 $|$      1.9343+1} &\multirow{2}{*}{16 $|$ 36:32}
\\[2pt] 
\mc{1}{|c|}{1030;6435} & & & & & \\[2pt] 
\hline\mc{1}{|c|}{dee10} 
	 &\multirow{2}{*}{0.068}	 &\multirow{2}{*}{4}&\multirow{2}{*}{ 8.4-7 $|$  \emph{2.9-1}} &\multirow{2}{*}{ 2.7107-1 $|$      3.3082-1} &\multirow{2}{*}{10 $|$ 17:04}
\\[2pt] 
\mc{1}{|c|}{365;8008} & & & & & \\[2pt] 
\hline\mc{1}{|c|}{friedman10} 
	 &\multirow{2}{*}{0.019}	 &\multirow{2}{*}{10}&\multirow{2}{*}{ 6.5-7 $|$  \emph{6.6-1}} &\multirow{2}{*}{ 5.9788-1 $|$      2.0224+0} &\multirow{2}{*}{14 $|$ 20:32}
\\[2pt] 
\mc{1}{|c|}{1200;3003} & & & & & \\[2pt] 
\hline\mc{1}{|c|}{mortgage5} 
	 &\multirow{2}{*}{0.001}	 &\multirow{2}{*}{70}&\multirow{2}{*}{ 5.3-7 $|$  \emph{6.2-1}} &\multirow{2}{*}{ 3.1051-2 $|$      4.4088-2} &\multirow{2}{*}{5:52 $|$ 1:21:45}
\\[2pt] 
\mc{1}{|c|}{1049;15504} & & & & & \\[2pt] 
\hline\mc{1}{|c|}{puma32h51tst3} 
	 &\multirow{2}{*}{0.012}	 &\multirow{2}{*}{8}&\multirow{2}{*}{ 1.9-7 $|$  \emph{6.5-2}} &\multirow{2}{*}{ 4.7014-3 $|$      1.1422-2} &\multirow{2}{*}{17 $|$ 59:00}
\\[2pt] 
\mc{1}{|c|}{1639;6545} & & & & & \\[2pt] 
\hline\mc{1}{|c|}{wizmir7} 
	 &\multirow{2}{*}{0.043}	 &\multirow{2}{*}{3}&\multirow{2}{*}{ 2.6-7 $|$  \emph{4.8-1}} &\multirow{2}{*}{ 7.1140-1 $|$      6.0186+0} &\multirow{2}{*}{22 $|$ 1:26:45}
\\[2pt] 
\mc{1}{|c|}{1461;11440} & & & & & \\[2pt] 
\hline

\end{tabular}
\end{center}
\end{table}

The second nonconvex problem we tested is the tuning-free robust regression problem with the MCP regularization. The performances are listed in Tables \ref{table-compare-random-MCP} and \ref{table-compare-KEEL-MCP}. For the synthetic datasets, we can see that the MCP regularizer can recover the true nonzero positions for the first 8 examples. Comparing with the $\ell_{1}$ regularizer, the MCP regularizer also achieves a sparse estimation with a fewer number of nonzeros for the synthetic type of data set. The corresponding solution is also much closer to $\dot{\beta}$ with less values of $L_{1}$, $L_{2}$, ME and FP. The PPMM algorithm can solve all these problems with high accuracy. In some results, the objective value obtained by the ADMM is smaller than that of the PPMM, but the solutions derived by the PPMM algorithm have fewer nonzeros. For example, the nnz for the ADMM are 7845, 3613, 11594, 169, 1224, 42 for the datasets example 12, baseball5, concrete7, dee10, friedman10 and mortgage5, respectively. The reason may be that Stage 1 of our PPMM algorithm generates a relative better sparse initial point.
\newpage
\begin{table}[!htbp]
\scriptsize
\begin{center}
\setlength{\belowcaptionskip}{10pt}
\parbox{.7\textwidth}{\caption{The performances of the PPMM and ADMM on the synthetic datasets for the MCP regularization. In this table,``$a$"=PPMM, ``$b$''=ADMM.}\label{table-compare-random-MCP}}
\begin{tabular}{|p{0.8cm}<{\centering}|p{0.5cm}<{\centering}|p{0.3cm}<{\centering}|p{1.8cm}<{\centering}|p{2.8cm}<{\centering}|p{1.5cm}<{\centering}|p{0.55cm}<{\centering}|p{0.55cm}<{\centering}|p{0.55cm}<{\centering}|p{0.3cm}<{\centering}|p{0.3cm}<{\centering}|}
\hline
\multirow{2}{*}{pbname}&\multirow{2}{*}{$\lambda$}&\multirow{2}{*}{nnz}&\multicolumn{1}{c|}{$\eta_{kkt}$}&\multicolumn{1}{c|}{pobj}&\multicolumn{1}{c|}{time}&\multirow{2}{*}{$L_{1}$}&\multirow{2}{*}{$L_{2}$}&\multirow{2}{*}{ME}&\multirow{2}{*}{FP}&\multirow{2}{*}{FN}\\\cline{4-6} $n$; $p$ & & & $a\ |\ b $ & $a\ |\ b $ & $a\ |\ b $ & & & & &\\[2pt] 
\hline\mc{1}{|c|}{exmp1} 
	 &\multirow{2}{*}{0.386}	 &\multirow{2}{*}{3}&\multirow{2}{*}{ 3.3-7 $|$  1.0-6} &\multirow{2}{*}{ 5.5706-1 $|$ 5.5717-1} &\multirow{2}{*}{41 $|$ 37:13} & \multirow{2}{*}{3.4-2} & \multirow{2}{*}{2.0-2} & \multirow{2}{*}{2.1-4} & \multirow{2}{*}{0} &\multirow{2}{*}{0}
\\[2pt] 
\mc{1}{|c|}{2000;8000} & & & & & & & & & &\\[2pt] 
\hline\mc{1}{|c|}{exmp2} 
	 &\multirow{2}{*}{0.158}	 &\multirow{2}{*}{3}&\multirow{2}{*}{ 8.9-7 $|$  \emph{2.3-3}} &\multirow{2}{*}{ 6.3745-1 $|$ 7.2690-1} &\multirow{2}{*}{26 $|$ 1:03:54} & \multirow{2}{*}{5.6-2} & \multirow{2}{*}{3.5-2} & \multirow{2}{*}{2.3-3} & \multirow{2}{*}{0} &\multirow{2}{*}{0}
\\[2pt] 
\mc{1}{|c|}{2000;8000} & & & & & & & & & &\\[2pt] 
\hline\mc{1}{|c|}{exmp3} 
	 &\multirow{2}{*}{0.200}	 &\multirow{2}{*}{3}&\multirow{2}{*}{ 8.7-7 $|$  \emph{4.6-6}} &\multirow{2}{*}{ 1.2375+0 $|$ 1.2377+0} &\multirow{2}{*}{19 $|$ 1:03:48} & \multirow{2}{*}{2.2-1} & \multirow{2}{*}{1.3-1} & \multirow{2}{*}{1.0-2} & \multirow{2}{*}{0} &\multirow{2}{*}{0}
\\[2pt] 
\mc{1}{|c|}{2000;8000} & & & & & & & & & &\\[2pt] 
\hline\mc{1}{|c|}{exmp4} 
	 &\multirow{2}{*}{0.101}	 &\multirow{2}{*}{3}&\multirow{2}{*}{ 9.2-7 $|$  \emph{1.9-1}} &\multirow{2}{*}{ 2.9299+0 $|$ 2.3082+1} &\multirow{2}{*}{16 $|$ 1:03:57} & \multirow{2}{*}{1.3-1} & \multirow{2}{*}{1.2-1} & \multirow{2}{*}{1.3-2} & \multirow{2}{*}{0} &\multirow{2}{*}{0}
\\[2pt] 
\mc{1}{|c|}{2000;8000} & & & & & & & & & &\\[2pt] 
\hline\mc{1}{|c|}{exmp5} 
	 &\multirow{2}{*}{0.157}	 &\multirow{2}{*}{3}&\multirow{2}{*}{ 3.6-7 $|$  \emph{8.7-6}} &\multirow{2}{*}{ 1.0466+0 $|$ 1.0470+0} &\multirow{2}{*}{26 $|$ 1:04:06} & \multirow{2}{*}{8.9-2} & \multirow{2}{*}{5.3-2} & \multirow{2}{*}{1.5-3} & \multirow{2}{*}{0} &\multirow{2}{*}{0}
\\[2pt] 
\mc{1}{|c|}{2000;8000} & & & & & & & & & &\\[2pt] 
\hline\mc{1}{|c|}{exmp6} 
	 &\multirow{2}{*}{0.183}	 &\multirow{2}{*}{3}&\multirow{2}{*}{ 6.6-7 $|$  \emph{5.9-6}} &\multirow{2}{*}{ 7.3170+0 $|$ 1.1263+2} &\multirow{2}{*}{22 $|$ 1:02:04} & \multirow{2}{*}{1.6-1} & \multirow{2}{*}{1.2-1} & \multirow{2}{*}{9.1-3} & \multirow{2}{*}{0} &\multirow{2}{*}{0}
\\[2pt] 
\mc{1}{|c|}{2000;8000} & & & & & & & & & &\\[2pt] 
\hline\mc{1}{|c|}{exmp7} 
	 &\multirow{2}{*}{0.069}	 &\multirow{2}{*}{25}&\multirow{2}{*}{ 3.0-7 $|$  \emph{4.5-2}} &\multirow{2}{*}{ 2.4807-1 $|$ 1.8024+0} &\multirow{2}{*}{51 $|$ 1:03:58} & \multirow{2}{*}{1.5-1} & \multirow{2}{*}{3.6-2} & \multirow{2}{*}{6.3-4} & \multirow{2}{*}{0} &\multirow{2}{*}{0}
\\[2pt] 
\mc{1}{|c|}{2000;8000} & & & & & & & & & &\\[2pt] 
\hline\mc{1}{|c|}{exmp8} 
	 &\multirow{2}{*}{0.054}	 &\multirow{2}{*}{25}&\multirow{2}{*}{ 1.1-7 $|$  \emph{1.6-1}} &\multirow{2}{*}{ 6.2953-1 $|$ 5.3778+0} &\multirow{2}{*}{41 $|$ 1:03:51} & \multirow{2}{*}{6.4-1} & \multirow{2}{*}{1.6-1} & \multirow{2}{*}{1.3-2} & \multirow{2}{*}{0} &\multirow{2}{*}{0}
\\[2pt] 
\mc{1}{|c|}{2000;8000} & & & & & & & & & &\\[2pt] 
\hline\mc{1}{|c|}{exmp9} 
	 &\multirow{2}{*}{0.036}	 &\multirow{2}{*}{25}&\multirow{2}{*}{ 2.6-7 $|$  \emph{9.6-2}} &\multirow{2}{*}{ 1.1546+0 $|$ 2.9945+0} &\multirow{2}{*}{32 $|$ 1:02:37} & \multirow{2}{*}{1.8+0} & \multirow{2}{*}{4.8-1} & \multirow{2}{*}{1.1-1} & \multirow{2}{*}{1} &\multirow{2}{*}{1}
\\[2pt] 
\mc{1}{|c|}{2000;8000} & & & & & & & & & &\\[2pt] 
\hline\mc{1}{|c|}{exmp10} 
	 &\multirow{2}{*}{0.046}	 &\multirow{2}{*}{21}&\multirow{2}{*}{ 5.8-7 $|$  \emph{6.3-2}} &\multirow{2}{*}{ 2.9500+0 $|$ 4.7911+0} &\multirow{2}{*}{26 $|$ 1:03:38} & \multirow{2}{*}{5.7+0} & \multirow{2}{*}{1.5+0} & \multirow{2}{*}{1.1+0} & \multirow{2}{*}{1} &\multirow{2}{*}{5}
\\[2pt] 
\mc{1}{|c|}{2000;8000} & & & & & & & & & &\\[2pt] 
\hline\mc{1}{|c|}{exmp11} 
	 &\multirow{2}{*}{0.039}	 &\multirow{2}{*}{26}&\multirow{2}{*}{ 7.5-7 $|$  \emph{1.2-1}} &\multirow{2}{*}{ 1.0282+0 $|$ 3.5867+0} &\multirow{2}{*}{46 $|$ 1:01:48} & \multirow{2}{*}{1.4+0} & \multirow{2}{*}{3.5-1} & \multirow{2}{*}{5.6-2} & \multirow{2}{*}{1} &\multirow{2}{*}{0}
\\[2pt] 
\mc{1}{|c|}{2000;8000} & & & & & & & & & &\\[2pt] 
\hline\mc{1}{|c|}{exmp12} 
	 &\multirow{2}{*}{0.031}	 &\multirow{2}{*}{41}&\multirow{2}{*}{ 1.8-8 $|$  \emph{2.3-5}} &\multirow{2}{*}{ 1.1617+1 $|$ 6.9967+0} &\multirow{2}{*}{2:12 $|$ 1:11:53} & \multirow{2}{*}{5.3+0} & \multirow{2}{*}{9.8-1} & \multirow{2}{*}{4.9-1} & \multirow{2}{*}{16} &\multirow{2}{*}{0}
\\[2pt] 
\mc{1}{|c|}{2000;8000} & & & & & & & & & &\\[2pt] 
\hline

\end{tabular}
\end{center}
\end{table}

\begin{table}[!htbp]
\begin{center}
\scriptsize
\setlength{\belowcaptionskip}{10pt}
\parbox{.7\textwidth}{\caption{The performances of the PPMM and ADMM on the KEEL datasets for the MCP regularization. In the table,``$a$"=PPMM, ``$b$''=ADMM.}\label{table-compare-KEEL-MCP}}
\begin{tabular}{|c|c|c|c|c|c|c|}
\hline
\multirow{2}{*}{probname}&\multirow{2}{*}{$\lambda$}&\multirow{2}{*}{nnz}&\multicolumn{1}{c|}{$\eta_{kkt}$}&\multicolumn{1}{c|}{pobj}&\multicolumn{1}{c|}{time}\\\cline{4-6} $n$; $p$ & & & $a\ |\ b$ & $a\ |\ b$ & $a\ |\ b$\\[2pt] 
\hline\mc{1}{|c|}{baseball5} 
	 &\multirow{2}{*}{0.071}	 &\multirow{2}{*}{30}&\multirow{2}{*}{ 9.0-7 $|$  \emph{3.0-5}} &\multirow{2}{*}{ 3.6849+2 $|$      9.6684+1} &\multirow{2}{*}{07 $|$ 36:03}
\\[2pt] 
\mc{1}{|c|}{337;20349} & & & & & \\[2pt] 
\hline\mc{1}{|c|}{compactiv51tra4} 
	 &\multirow{2}{*}{0.061}	 &\multirow{2}{*}{3}&\multirow{2}{*}{ 6.0-7 $|$  \emph{2.9-3}} &\multirow{2}{*}{ 4.6827+0 $|$      1.3629+1} &\multirow{2}{*}{3:58 $|$ 4:01:07}
\\[2pt] 
\mc{1}{|c|}{6553;12650} & & & & & \\[2pt] 
\hline\mc{1}{|c|}{concrete7} 
	 &\multirow{2}{*}{0.011}	 &\multirow{2}{*}{24}&\multirow{2}{*}{ 7.5-7 $|$  \emph{7.5-6}} &\multirow{2}{*}{ 3.1083+0 $|$      1.5686+0} &\multirow{2}{*}{15 $|$ 35:01}
\\[2pt] 
\mc{1}{|c|}{1030;6435} & & & & & \\[2pt] 
\hline\mc{1}{|c|}{dee10} 
	 &\multirow{2}{*}{0.020}	 &\multirow{2}{*}{7}&\multirow{2}{*}{ 4.6-7 $|$  \emph{1.5-6}} &\multirow{2}{*}{ 2.1662-1 $|$      1.8607-1} &\multirow{2}{*}{06 $|$ 17:03}
\\[2pt] 
\mc{1}{|c|}{365;8008} & & & & & \\[2pt] 
\hline\mc{1}{|c|}{friedman10} 
	 &\multirow{2}{*}{0.010}	 &\multirow{2}{*}{12}&\multirow{2}{*}{ 4.9-7 $|$  \emph{6.5-6}} &\multirow{2}{*}{ 5.8870-1 $|$      3.7554-1} &\multirow{2}{*}{07 $|$ 20:27}
\\[2pt] 
\mc{1}{|c|}{1200;3003} & & & & & \\[2pt] 
\hline\mc{1}{|c|}{mortgage5} 
	 &\multirow{2}{*}{0.053}	 &\multirow{2}{*}{6}&\multirow{2}{*}{ 3.0-7 $|$  1.0-6} &\multirow{2}{*}{ 1.6327-1 $|$      1.5034-1} &\multirow{2}{*}{32 $|$ 39:48}
\\[2pt] 
\mc{1}{|c|}{1049;15504} & & & & & \\[2pt] 
\hline\mc{1}{|c|}{puma32h51tst3} 
	 &\multirow{2}{*}{0.015}	 &\multirow{2}{*}{3}&\multirow{2}{*}{ 7.8-7 $|$  \emph{1.4-2}} &\multirow{2}{*}{ 4.3927-3 $|$      1.9269-2} &\multirow{2}{*}{17 $|$ 59:08}
\\[2pt] 
\mc{1}{|c|}{1639;6545} & & & & & \\[2pt] 
\hline\mc{1}{|c|}{wizmir7} 
	 &\multirow{2}{*}{0.046}	 &\multirow{2}{*}{6}&\multirow{2}{*}{ 3.4-7 $|$  \emph{3.2-6}} &\multirow{2}{*}{ 7.1795-1 $|$      2.6756+0} &\multirow{2}{*}{25 $|$ 1:25:06}
\\[2pt] 
\mc{1}{|c|}{1461;11440} & & & & & \\[2pt] 
\hline

\end{tabular}
\end{center}
\end{table}

\newpage
\section{Conclusion}
\label{sec:Conclusion}
In this paper, we introduce the PPMM algorithm for solving nonconvex tuning-free robust regression problems. We obtain the solution by solving a sequence of convex majorization-minimization problems. The SSN based PPA is applied to solve each convex subproblem. We make full use of the structure of the Clarke generalized Jacobian and solve each subproblem efficiently. We also prove that the PPMM algorithm converges to a d-stationary point and analyze the convergence rate based on the KL property of the problem.


\end{document}